 \documentclass[draft]{article}

\usepackage{amsmath,amsfonts,amsthm,amssymb,amscd,cancel,color,mathabx}
\usepackage{enumitem}
\usepackage{verbatim}
\usepackage[dvipdfmx]{graphicx}
\usepackage{ulem,color}

\renewcommand{\Re}{\mathop{\rm Re}\nolimits}
\renewcommand{\Im}{\mathop{\rm Im}\nolimits}
\def\S{\mathhexbox278}

\theoremstyle{plain}
\newtheorem{theorem}{Theorem}[section]
\newtheorem{lemma}[theorem]{Lemma}
\newtheorem{proposition}[theorem]{Proposition}

\theoremstyle{definition}

\theoremstyle{remark}
\newtheorem{remark}[theorem]{Remark}
\newtheorem{notation}[theorem]{Notation}

\newtheorem{claim}[theorem]{Claim}
\newtheorem{statement}[theorem]{Statement}

\newcommand{\R}{{\mathbb R}}

\def\im{{\rm i}}

\newcommand{\C}{\mathbb{C}}

\def\({\left(}
\def\){\right)}
\def\<{\left\langle}
\def\>{\right\rangle}
\newcommand{\sech}{{\mathrm{sech}}}

\newcommand{\supp}{{\mathrm{supp}\ }}
\newcommand{\rad}{{\mathrm{rad}}}    \newcommand{\diag}{{\mathrm{diag}}}
\newcommand{\Span}{{\mathrm{Span}}}

\newtheorem{case}{Case}
\makeatletter

\numberwithin{equation}{section}

\begin{document}

\title{On the  asymptotic stability  on the line  of  ground states of the  pure power NLS with  $   0\le 2-p \ll 1 $.}

\author{Scipio Cuccagna, Masaya Maeda }
\maketitle

\begin{abstract}We continue our series  devoted, after references \cite{CM24D1} and \cite{CM243}, at proving the asymptotic stability of ground states of the pure power Nonlinear Schr\"odinger equation on the line. Here we assume  some results on the spectrum of the linearization obtained computationally by Chang et al. \cite{Chang} and then we explore the equation for exponents $p\le 2$ sufficiently close to 2. The ensuing loss of regularity of the nonlinearity  requires new arguments.
\end{abstract}

\section{Introduction}

We consider  the pure    power focusing  Nonlinear Schr\"odinger Equation  (NLS) on the line
\begin{align}\label{eq:nls1}&
  \im \partial _t  u +\partial _x^2  u = -f(u)  \text{   where }    f(u)=|u| ^{p-1}u    \text{ for $ 0\le 2-p \ll 1$. }
\end{align}
      It is well known that equation   \eqref{eq:nls1} has solitary waves
     \begin{equation}\label{eq:solwave}
      u(t,x)= e^{\frac \im 2 v  x
-\frac \im 4 |v|^2 t +\im t\omega +\im \gamma } \phi _\omega
(x-vt- D)
     \end{equation}
         where  $\phi _\omega (x)=\omega ^{\frac 1{p-1}} \phi (\sqrt{\omega }x) $  with, formula (3.1) of Chang et al. \cite{Chang},
\begin{align} &
    \phi (x) =  {\(\frac {p+1}2 \)^{\frac 1{p-1}}}{
\sech   ^{\frac 2{p-1}}\(\frac{p-1}2 x\)}  .\label{eq:sol}
   \end{align}
Energy $\mathbf{E}$,  mass $\mathbf{Q}$ and linear momentum $\mathbf{P} $  are invariants of \eqref{eq:nls1}, where
  \begin{align}\label{eq:energy}
& \mathbf{E}( {u})=\frac{1}{2}   \| u' \| ^2 _{L^2\( \R \) } -\int_{\R} F(u)    \,dx   \text {  where }  F(u)    =\dfrac{ |u|  ^{p+1}}{p+1}, \\&   \label{eq:mass} \mathbf{Q}( {u})=\frac{1}{2}   \| u  \| ^2 _{L^2\( \R \) }  \text{   and } \\&  \mathbf{P}  ( {u})  =-\frac{1}{2}  \< \im \partial _x u, u \>   \  \text{ where }  \< f, g \>  =\Re \int _\R f(x) \overline{g}(x) dx.\label{eq:moment}
\end{align}
 It is well known that
   $\phi _\omega$  minimizes $\mathbf{E}$ under the constraint $\mathbf{Q}=\mathbf{Q}(\phi _\omega )=:\mathbf{q}(\omega )$. Notice that $\mathbf{q}(\omega ) =\omega ^{ \frac{2}{p-1}-\frac{1}{2}}   \mathbf{q}(1)$.
In particular, we have $\nabla \mathbf{E} (\phi _\omega)= -\omega   \nabla \mathbf{Q} (\phi _\omega)  $ which reads also
\begin{align} \label{eq:static}
   - \phi _\omega ''+\omega \phi _\omega - \phi _\omega ^p=0   .
\end{align}
Consider now  for $\omega ,\delta \in \R_+:=(0,\infty)$ the set     $$\mathcal{U} (\omega  ,\delta  ) := \bigcup _{\vartheta _0 , x_0\in \R }   e^{\im \vartheta _0} e^{-x_0 \partial _x}  D_{H ^  (\R )}({\phi}_{\omega },\delta  ), $$ where $D_X(u,r):=\{v\in X\ |\ \|u-v\|_X<r \}$.
The following was shown by Cazenave and Lions \cite{cazli}, see also Shatah  \cite{shatah} and Weinstein  \cite{W1}.

\begin{theorem}[Orbital Stability]  Let $p\in (1,5)$ and let $\omega _0 >0$. Then for any $\epsilon >0$  there exists a  $\delta >0$ such that for any  initial value    $u_0\in \mathcal{U} (\omega _0,\delta  ) $ then   the corresponding solution  satisfies
$u\in C^0  \(  \R ,   \mathcal{U} (\omega _0,\epsilon  ) \) $.
\end{theorem}
\qed

Notice that, as we discuss  after formula \eqref{eq:orbsta1}, it is possible to get $\delta\sim \epsilon $, as discussed in Stuart \cite{stuart}, but this observation is not essential here.
In order to study the notion of asymptotic stability, like in finite dimension, it is useful to have information on the  \textit{linearization} of \eqref{eq:nls1} at $\phi _\omega$, which we will see later has the following form
\begin{align}\label{eq:lineariz2} \partial _t   \begin{pmatrix}
 r_1 \\ r_2
 \end{pmatrix}      &=  \mathcal{L}_{\omega }  \begin{pmatrix}
 r_1 \\ r_2
 \end{pmatrix}  \text{  with }    \mathcal{L}_\omega := \begin{pmatrix}
0 & L_{-\omega} \\ -L_{+\omega} & 0
\end{pmatrix} ,
    \end{align}
where    \begin{align}\label{eq:lin1}&
  L _{+\omega}:=- \partial _x^2   +\omega -p \phi _\omega ^{p-1}   \\& \label{eq:lin0} L _{-\omega}:=- \partial _x^2   +\omega-  \phi _\omega ^{p-1}.
\end{align}
The linearization is better seen in the context of functions  in $H^ 1 (\R , \R ^2 ) $ rather than in $H^ 1 (\R ,\C ) $, because it is $\R$--linear rather than $\C$--linear.    Defining the scaling
\begin{align}
  \label{scaling}  S_{\omega}u(x)=\omega^{1/4}u(\sqrt{\omega}x) \Longrightarrow   \mathcal{L}_{\omega} = \omega S_{\omega} \mathcal{L} S_{\omega}^{-1} \text{ for } \mathcal{L}=\mathcal{L}_1.
\end{align}
It is well known that the spectrum of $\mathcal{L} $ is symmetric with respect to real and imaginary axes and that the essential spectrum is $(-\im \infty,-\im]\cup[\im ,\im \infty)$.
 Chang et al. \cite{Chang} compute  numerically the   statements   \textbf{{1}}-- \textbf{{5}}.  Coles and Gustafson \cite{coles}   prove  rigorously \textbf{{6}}. Krieger and Schlag \cite{KrSch}  prove rigorously \textbf{{7}}.

 \begin{statement} \label{stat:eigen}  \begin{description}
                      \item[1]  There is a value  $p_0$ close to $1.82$ such that  for $p_0<p<3$ and for  $3<p<5$ the linearization  $\mathcal{L} $ has exactly one eigenvalue of the form $\im \lambda(p)$  with $0<\lambda(p)<1$.
                       \item[2]   For $p_0<p<3$    we have  $ \lambda(p)>1/2$.
                       \item[3] There is a $p_1\in \( 1, p_0 \)$ such that  for $p_1 <p< p_0$  the linearization   $\mathcal{L} $ has  a second eigenvalue
 $\im \mu (p)$ appears with  $0<\lambda(p)<\mu (p) <1$.
 \item[4]  The function $(3,5)\ni p\mapsto \lambda (p) $   is  strictly decreasing with $\lambda (3^ {+})=1$, $\lambda(4)>1/2$  and $\lambda (5^{-} )=0 $.

      \item[5]  The function $(1,3)\ni p\mapsto \lambda (p) $   is  strictly increasing with $\lambda (3^ {-})=1$  and $\lambda (1^+ )=0 $.

      \item[6]   For $0<|p-3|\ll 1$      the linearization   $\mathcal{L} $ has       exactly one eigenvalue of the form $\im \lambda(p)$  with
       \begin{equation}\label{eq:eigpnear3}
         \lambda(p) \sim 1 - \alpha _0 (p-3) ^4 +o\( (p-3) ^4  \) \text{ for a constant $\alpha _0>0$.} \end{equation}

          \item[7]  There are no eigenvalues of the form $\im \mu$ with $\mu >1$.

                    \end{description}

 \end{statement}

We  have the following,   which we proved in \cite{CM243}.
\begin{proposition}[No Threshold Resonance]
   \label{prop:nonres}   Consider the following statement,
\begin{equation}\label{eq:nores111}
  \text{ there is   no nonzero bounded  solution of $\mathcal{L} u=\im   u$.}
\end{equation}
Consider the set
\begin{equation}\label{eq:deftildeF}
  \widetilde{{\mathbf{F}}} :=\{p \in ( 1, +\infty ) \ |\  \text{statement \eqref{eq:nores111} is true }\} .
\end{equation}
Then  $( 1,  +\infty ) \backslash   \widetilde{{\mathbf{F}}}$ is a discrete subset of $( 1, +\infty )$.
\end{proposition} \qed

An important ingredient in this paper is the  nonlinear Fermi Golden Rule(FGR). Recall the  $p_0$   in Statement \ref{stat:eigen}.
 \begin{proposition} \label{prop:FGR} Consider the constant $\gamma (p) $  introduced later in \eqref{eq:fgrgamma}. Then assuming as hypotheses \textbf{{1}}-- \textbf{{2}}
 in Statement \ref{stat:eigen} we have that $\gamma (\cdot )\in C ^{\omega} ((p_0,3) , \R )$ and that the set
    \begin{equation}\label{eq:deftildeF}
   {{\mathbf{F}}} :=\{p \in   (p_0,3) \ |\ \gamma (p)   \neq 0 \}
\end{equation} is s.t. $(p_0,3) \backslash  {{\mathbf{F}}}$ is a discrete subset of $(p_0,3)$.
 \end{proposition}
\proof The fact that $\gamma (\cdot )\in C ^{\omega} ((p_0,3) , \R )$, i.e. $\gamma (p)$ is analytic in $p$, follows from the fact that $\mathcal{L}$ is analytic of type $(A)$  in $p$, as we will see later.
 Then the last line of the statement follows from the fact, proved in  \cite{CM24D1,CM242}, that it is possible to define  $\gamma (p)$ so that it is nonzero for   $0<|p-3|\ll 1$.  \qed

In this paper we prove  the following result.

\begin{theorem} \label{thm:asstab} There is a   $\widetilde{p}_0\in (p_0, 2)$  such that if $ {p}\in (\widetilde{p}_0, 2]$ belongs to  $\widetilde{{\mathbf{F}}} \cap  {{\mathbf{F}}}$ then for any $\omega _0 >0$, any $a>0$ and any $\epsilon >0$   there exists a  $\delta >0$ such that for any  initial value    $u_0\in  D _{H^1 (\R ) } ({\phi}_{\omega  _0},\delta  ) $
there exist functions $ \vartheta , \omega ,v, D \in C^1 \( \R , \R \)$,  $ z \in C^1 \( \R , \C  \)$,  $\omega_+>0$  and $v_+\in \R$    s.t.\  the  solution  of \eqref{eq:nls1} with initial datum  $u_0$  can be written as
 \begin{align} \label{eq:asstab1}
    & u(t)= e^{\im \vartheta (t)-D(t) \partial _x} \(   \phi _{\omega (t)} + z(t) \xi _{\omega(t) }+ \overline{z}(t) \overline{\xi} _{\omega (t)} +\eta (t)\) \text{   with}
 \\&  \label{eq:asstab2}   \int _{\R } \|  e^{- a\< x\>}   \eta (t ) \| _{H^1(\R )}^2 dt <  \epsilon  \text{  where }\< x\>:=\sqrt{1+x^2}, \\&  \label{eq:asstab20}
   \lim _{t\to +\infty  }  \|  e^{- a\< x\>}   \eta (t ) \| _{L^2(\R )}     =0  , \\&  \label{eq:asstab2}
   \lim _{t\to +\infty}z(t)=0      ,\\&
\lim _{t\to +\infty}\omega (t)= \omega _+  \text { and} \label{eq:asstab3} \\& \lim _{t\to +\infty}v (t)= v _+  . \label{eq:asstab4}
\end{align}
Here $\widetilde{p}_0$, besides what stated  above is also such that
\begin{align}
  \label{eq:tildep0}  2-p - \sqrt{1-\lambda (p)}< 0\text{ for all $ {p}\in (\widetilde{p}_0, 2]$.}
\end{align}
\end{theorem}

\begin{remark}\label{rem:asstab3--}  Notice that the fact that the $H^1(\R )$ norm of $\eta$ is uniformly bounded for all times, guaranteed by the orbital stability, and
Theorem \ref{thm:asstab}  imply   that
\begin{align}\label{eq:rem:asstab3--1}
  \lim _{t\to +\infty }  u(t) e^{-\im \vartheta (t)+ D(t) \partial _x  }=  \phi  _{\omega _+}  \text{  in }L^\infty _{\text{loc}}(\R ) .
\end{align}
\end{remark}

Equation \eqref{eq:nls1} is a classical Hamiltonian systems in PDE's and the asymptotic stability of its ground states has been a longstanding open problem considered the early papers  by Soffer  and Weinstein \cite{SW1,SW2}.
For a rather extensive discussion of the problem, we refer to \cite{CM24D1}, where we deal with the  case $p\sim 3$ but $p\neq 3$ (for $p=3$  Theorem \ref{thm:asstab} is false).
  The results in \cite{CM24D1} are only partially   based on the condition that $p$ is close to 3 and can be extended to a larger set of $p$'s thanks to an analogue of  Proposition  \ref{prop:nonres}   and an analogue of Proposition  \ref{prop:FGR}. However the proof in
\cite{CM24D1} exploits in various places the condition $p>2$, which guarantees  the $C^2 $ regularity of the function $u\to f(u)$.  The novelty of this paper is the fact that we tackle low regularity  cases with $p\le 2$,  utilizing ideas similar to  \cite{CMMS23}.  Notice that the key for the stabilization  is given by the  linear dispersion of the continuous mode  $\eta$ and by the friction exerted on the   discrete mode $z$ contained in \eqref{eq:asstab1} by the nonlinear interaction with the continuous mode  $\eta$ by means of the nonlinear Fermi Golden Rule (FGR) (here $z$ is like the variable of a finite dimensional  hamiltonian system which is embedded and not well isolated inside an infinite dimensional one, so it is slowly losing all its energy; the notion was introduced by Sigal \cite{sigal}). The fact that dispersion occurs even with such low values of $p$  is remarkable, in view of how difficult and central has been in the literature the problem of dealing with dispersion in the presence of a nonlinearity, which has lead to milestones like Klainerman's vector fields and   null forms, see   for instance \cite{Kl86},  Shatah's normal forms    \cite{shatahnf} and Christodoulou's conformal mapping   \cite{Christodoulou}, see also the discussion in Strauss \cite{strauss}. Notice that   a version of
Theorem   \ref{thm:asstab}, with $u=\eta$, is false near 0 for  $p < 5$  and a counterexample to \eqref{eq:asstab20} is given   by the ground states
$u(t,x) = e^{ \im t\omega   }   \phi  _{\omega}$ since $\|   \phi  _{\omega} \| _{H^1} \xrightarrow{\omega \to 0^+} 0 $. For $p\ge 3$ by Ozawa \cite{ozawa} for any $u_+\in L^{2,2}(\R )$, the space defined in  \eqref{eq:l2w} below, and both for focusing and defocusing equation   there exists a solution $u$ of \eqref{eq:nls1} with $\| u(t) - e^{\im S(t,\cdot )}e^{\im t\partial ^2_x}u_+\| _{L^2(\R )}\le O(t^{-\alpha}) $ with $\alpha >1/2$ and for an appropriate real phase $S(t,x )$, which is  equal to 0 if $p>3$. Since   for  $|x|\le 1$   we have $e^{\im t\partial ^2_x}u_+( x) = \frac{1}{ \sqrt{2\im t}  } e^{\im \frac{x^2}{4t}} \widehat{u}_+(0) +O(t^{-\alpha}) $ with $\alpha >1/2$, in the generic case with   $\widehat{u}_+(0)\neq 0$  we conclude that $\|   u  \| _{L^2 (\R _+ , L^2(|x|\le 1))} =+\infty$. For $p\in (1,3)$  it is known that
nonzero solutions of the defocusing equation behave differently from the solutions of the linear equation and we do not know if a result like   \eqref{eq:asstab20}  with $u=\eta$ is false.
The above considerations show   that generally
   equation \eqref{eq:nls1}  displays better dispersive properties near the ground state than near 0.  In fact, as we have seen,  there is more dispersion near the soliton than near vacuum for the constant coefficient linear equation, due to
   the fact that the linearization at the ground state is better behaved in terms of the Kato smoothing, that is \cite{kato}, in view of  {Proposition}
   \ref{prop:nonres}, while  operator $-\partial _x^2$     has a resonance.

    We state here briefly the main ingredients of the proof.  First of all, we use our notion of Refined Profile to find a  particularly favourable ansatz for our solution.  Thanks to this ansatz we obtain various convenient cancellations  and good estimates for the discrete components of our solution.  The ansatz is also particularly convenient to prove the FGR, both in the rather simple spectral configuration considered here, but also in much more complicated situations like for example in \cite{CMS23}.  Next we use the first virial inequaity   
Kowalczyk et al. \cite{KM22,KMM20,KMM3,KMMvdB21AnnPDE}, which  can traced to the work of Martel and Merle, see for example  \cite{MarMerle} and therein, and Merle and Raphael, see for example  \cite{MR4} and therein, and   earlier in Kato \cite{katovirial}.   The crux is represented by the (almost) positive commutator formula \eqref{eq:poscomm} and by the remarkable pure power inequality  \eqref{eqpurepower}.   It is thanks to the flexibility and robustness of this tool  that we can deal with low values of   $p$. Finally we prove a low energy dispersive estimate.  Here, like in \cite{CM24D1},  we modify the approach of
Kowalczyk et al.  \cite{KMM3}, Martel \cite{Martel1,Martel2}  and    Rialland  \cite{rialland, Rialland2}  by  replacing their  second virial inequality with a smoothing estimate.   The advantage   compared to the is that  we    can deal with equations which are not small perturbations of the cubic NLS, contrary to   Martel \cite{Martel1,Martel2} and Rialland  \cite{rialland, Rialland2}.
As we discussed at length in \cite{CM24D1}, all  is needed is a smoothing estimate    for the continuous mode $\eta $  multiplied by a cutoff. We can derive it    easily by multiplying the equation of $\eta$ by a cutoff. In this way the term $f(\eta )$ once multiplied by the cutoff is easy to bound. Compared to \cite{CM24D1}, we benefit by a better smoothing estimate proved in \cite{CM243}.

 Besides the papers involving the virial inequalities 
there is also a large literature which uses
  dispersive estimates for   problems which are related but  in general are different from the  one we consider here, for example and without trying to give a complete list
\cite{DM20,GPR18,GP20,GPZ,GoPu,GoPu1,naumkin2016}, \cite{LP0}--\cite{LS2023}.   In general these methods seem to require quite regular $u\to f(u)$, so they do not seem applicable here. An interesting mixing of  phase space analysis and virial inequalities  in Palacios and Pusateri \cite{PalPus24}.
Germain and Collot \cite{germain2} have recovered and partially generalized Martel \cite{Martel1} by using very different methods.  Li and  Luhrmann \cite{LiLu} for the cubic NLS  partially recover  the result obtained in  \cite{cupe2014} (using the the integrable structure)   with methods which do not exploit  explicitly     the integrable structure of the cubic NLS.
A  rather long list of references on the asymptotic stability of ground states of the NLS   up until 2020 is in  our  survey \cite{CM21DCDS}.

\section{Linearization}\label{sec:lin1}

We return to a discussion of the linearization \eqref{eq:lineariz2}. First of all we notice
\begin{align}
   \label{eq:lambdap}  \partial_{\omega}  \( \omega^{\frac{1}{p-1}}u(\sqrt{\omega}x) \) =   \omega^{\frac{1}{p-1}-1}  (\Lambda_p   u)(\sqrt{\omega}x) \text{ for } \Lambda_p=\frac{1}{2}x\partial_x + \frac{1}{p-1}
 \end{align}
so that in particular  $\Lambda_p u =\left.\partial_{\omega}\right|_{\omega=1} \omega^{\frac{1}{p-1}}u(\sqrt{\omega}x)$.
Weinstein \cite{W2}  showed that  for $ 1<p<5$ the   generalized kernel $ {N}_g(\mathcal{L} ):=\cup_{j=1}^\infty \mathrm{ker} \mathcal{L} ^j$ in $H^1  ( \R , \C^2 ) $   is
\begin{align}\label{eq:Ng}
 {N}_g(\mathcal{L} )=\mathrm{span}\left \{ \begin{pmatrix}
 0 \\ \phi
 \end{pmatrix} ,\ \begin{pmatrix}
\partial_x    \phi   \\ 0
 \end{pmatrix}    ,\ \begin{pmatrix}
 \Lambda _p   \phi  \\ 0
 \end{pmatrix}      , \   \begin{pmatrix}
 0 \\ x\phi
 \end{pmatrix}  \right \} .
\end{align}
    By symmetry reasons,   it known that the spectrum $\sigma \(  \mathcal{L} \) \subseteq \C $ is symmetric by reflection with respect of the coordinate axes. Furthermore,
by Krieger and Schlag   \cite[p. 909]{KrSch}  we know that  $\sigma \(  \mathcal{L} \) \subseteq \im \R  $. By standard Analytic  Fredholm theory   the essential spectrum is $(- \infty \im ,-  \im ] \cup
[ \im   , +\infty \im ) $. As we mentioned in Statement \ref{stat:eigen},   Krieger and Schlag   \cite[pp. 914--915]{KrSch} prove that  $\mathcal{L} $    has  no eigenvalues  $\im \mu $ with $\mu >1 $.

\noindent Let us consider the orthogonal  decomposition
\begin{align}
  \label{eq:dirsum} L^2 (\R , \C ^2) =  {N}_g(\mathcal{L} )\bigoplus {N}_g^\perp (\mathcal{L} ^*)
\end{align}  We have,  for $\lambda  = \lambda (p,\omega)$,   a further decomposition
\begin{align}
  \label{eq:dirsum1}&{N}_g^\perp (\mathcal{L} ^*) =   \ker (\mathcal{L} -\im \lambda (p) )\bigoplus \ker (\mathcal{L} +\im \lambda (p)) \bigoplus X_c   \text{ where } \\& X_c  : = \(   {N}_g(\mathcal{L} ^*) \bigoplus  \ker (\mathcal{L} ^*-\im \lambda (p) )\bigoplus \ker (\mathcal{L}^*+\im \lambda _{\omega } ) \) ^{\perp} . \label{eq:dirsum11}
\end{align}
We denote by $P_c$ the projection of $ L^2  (\R , \C ^2)$  onto $X_c (\omega )$ associated with the above decompositions.  We have the following.
\begin{lemma}\label{lem:simandmult}
  Consider the  $p_0$ in Statement \ref{stat:eigen}. Then for $p\in (p_0,3)$  the algebraic and dimension of the eigenvalue $\im \lambda (p)$ is one   and $  \ker \( H -\lambda (p) \) \subset H^1 _\rad \( \R , \C^2\)$, see Notation \ref{not:notation} below.
   \end{lemma}
\proof  First of all, Coles and Gustafson \cite{coles}  show that $\dim \ker \( H -\lambda (p) \) =1$  for     $p $ close to 3.   Next,  we show later in Lemma  \ref{eq:multeig}
that $\dim \ker \( H -\lambda (p) \) = N_g\( H -\lambda (p) \)$  if $\lambda (p) \in (0,1)$. Since  $H$ is analytic of type $(A)$  in $p$ the eigenvalue $\lambda (p)$ depends analytically on $p \in (p_0, 3)$, see \cite[Ch. 12]{reedsimon}.   Since $\lambda (p)$ remains an isolated eigenvalue  for $p \in (p_0, 3)$, then $\dim \ker \( H -\lambda (p) \)$ remains constant.

Turning to the parity of the elements of $\ker \( H -\lambda (p) \)$, for $p$ close to 3 they  are indeed  even,    as is shown in \cite[Formula 11.8]{CM24D1}. Then again  we can conclude that this is true for all
 $p\in (p_0,3)$.

 \qed

\noindent The space $L^2 (\R , \C ^2)$ and the action of $\mathcal{L}_{\omega }$ on it is obtained by first  identifying $L^2 (\R , \C  ) =L^2 ( \R ,  \R ^2  )  $ and then by extending   this action to  the completion of
$ L^2  (\R ,  \R ^2  ) \bigotimes _\R \C  $ which is identified with $L^2  (\R , \C ^2)$.  In $\C$ we consider the inner product,
\begin{align*}
   \< z, w \> _{\C} =\Re \{z\overline{w}    \} =z_1w_1+z_2w_2  \text{ where } a_1=\Re a, \quad a_2=\Im a   \text{ for }a=z,w.
\end{align*}
This obviously coincides with the inner product in $\R ^2$ and expands as the standard sesquilinear  $  \< X  , Y  \> _{\C ^2}  =  X ^ \intercal \overline{Y}$  (row column product, vectors here are columns)  form   in $\C ^2$.  The operator of multiplication by $\im $  in $C=\R^2$ extends into the linear operator
$J ^{-1}=-J$ where   \begin{align}
J=\begin{pmatrix}
0 & 1 \\ -1 & 0
\end{pmatrix}.\nonumber
\end{align}
For  $u,v\in L^2 (\R , \C ^2)$ we set $ \< u , v  \>    :=\int _\R \< u (x), v (x) \> _{\C ^2}  dx$. We have a natural symplectic form given by $\Omega :=\<  J ^{-1}\cdot  , \cdot   \>$ in both $ L^2(\R , \C ^2)$ and
$ L^2 (\R , \R ^2)=L^2  (\R , \C  )$, where equation \eqref{eq:nls1} is the Hamiltonian system  in  $  L^2  (\R , \C  )$     with Hamiltonian the energy $E$ in \eqref{eq:energy}.  As we mentioned we consider
a generator $\xi \in  \ker (\mathcal{L  }-\im \lambda (p) )$.  Then for the complex conjugate $\overline{\xi}  \in  \ker (\mathcal{L} +\im \lambda (p))$.
Notice   the well known and elementary $J \mathcal{{L}} =-\mathcal{{L}}  ^*J$  implies that  $ \ker (\mathcal{L} ^*+\im \lambda  (p )) = \mathrm{span}\left \{  J {\xi}  \right \}$ and $ \ker (\mathcal{L} ^*-\im \lambda (p))= \mathrm{span}\left \{  J \overline{{\xi}}  \right \}$. All the above is discussed also in \cite{CM24D1} where we recall   that we can  use the normalization
\begin{align}
  \label{eq:xinormlize} \Omega (   {\xi}  ,  {\xi}  ) =- \im
\end{align} and that furthermore we can assume
\begin{align}\label{eq:reimxi1} &  \xi  = \(   \xi _1 , \xi _2 \) ^\intercal  \text{  with }  \xi _1=\Re  \xi _1 \text{  and }  \xi _2=\im \Im  \xi _2,
\end{align}
so that \eqref{eq:xinormlize} becomes
\begin{align}
  \label{eq:xinormliz1} \int _{\R} \xi _1 \Im \xi _2 dx = 2 ^{-1}.
\end{align}
Notice also that $\xi$ is even in $x$.

\begin{lemma}
  \label{lem:deceig} There exists  $\widehat{p}_0\in (p_0, 2)$ such that for $p\in ( \widehat{p}_0, 2]$  we have the point--wise estimate
  \begin{align}\label{eq:deceig1}
    |\xi (x)|\le C_p e^{-|x|\sqrt{1-\lambda (p)}}.
  \end{align}
\end{lemma}
\proof  Introducing
\begin{align}\label{def:U}
	U = \begin{pmatrix}
		1 &
		1  \\
		\im &
		-\im      \end{pmatrix} \, , \quad
	U^{-1}= \frac 12  \begin{pmatrix}  1 &
		-\im   \\
		1 &
		\im    \end{pmatrix}   \end{align}
	we have
$
 U^{-1} J
	U= \im \sigma _3  $
 so that
\begin{align}\label{eq:opH}
   & U ^{-1}\mathcal{L} _{\omega} U = \im  H _{\omega}   \text{  where } H _{\omega}  =   \sigma _3 \(-\partial _x^2 +\omega \) +V  _{\omega}, \\&
    V _{\omega} :=  \omega  M_0 \sech ^2   \(  \sqrt{\omega}  \frac{p-1}2 x\)  \text{ with } M_0= -\( \frac{p+1}{2}
    \sigma _3      +\im  \frac{p-1}{2}  \sigma _2 \)   \frac {p+1}2     \nonumber
\end{align}
where we recall the Pauli matrices \begin{equation}\label{eq:pauli}
   \sigma_1=\begin{pmatrix} 0 &
   1 \\
   1 & 0
    \end{pmatrix} \,,
   \quad
   \sigma_2= \begin{pmatrix} 0 &
   -\im \\
    \im  & 0
    \end{pmatrix} \,,
   \quad
   \sigma_3=\begin{pmatrix} 1 &
   0 \\
   0 & -1
    \end{pmatrix}.
   \end{equation}
Then for $H_0= \sigma _3 \(-\partial _x^2 +1 \) $,  $H= \left . H _{\omega} \right | _{\omega =1} $   and  $V= \left . V _{\omega} \right| _{\omega =1} $    we have
\begin{align*}
 \( H _0  -\lambda (p) + V   \)  w= 0 \text{ for }w:=U^{-1}\xi .
\end{align*}
It follows that, for    $\diag (a, b)$   the diagonal matrix with first $a$ and then $b$ on the diagonal,
\begin{align*}
  w = - \diag \(  R _{-\partial _x ^2}  ( \lambda (p)-1   )  ,    - R _{-\partial _x ^2}  ( -\lambda (p)-1   )    \) Vw
\end{align*}
so that pointwise we have
\begin{align*}
   &|w(x)|\lesssim \int _{\R} e^{-|x  -y|\sqrt{1-\lambda (p)}}  e^{-(p-1) |y|}|w(y)| dy   .
\end{align*}
Since $w\in L^\infty ( \R )$  by $\sqrt{1-\lambda (p)}  < 2^{-\frac{1}{2}} $  for $p>p_0$   we can pick $\widehat{p}_0$ so that  $\sqrt{1-\lambda (p)}  < p-1 $   for  $p\in (\widehat{p}_0, 2]$.
Then by the triangular inequality
\begin{align*}
   &  \int _{\R} e^{-|x  -y|\sqrt{1-\lambda (p)}}  e^{-(p-1) |y|}|w(y)| dy    \le   e^{-|x  |\sqrt{1-\lambda (p)}} \int _{\R}    e^{ \( \sqrt{1-\lambda (p)}  -(p-1)\) |y|}  dy \| w \| _{L^\infty}
\end{align*}
yields  the following and so also \eqref{eq:deceig1},
\begin{equation*}
     |w(x)|\lesssim   e^{-|x  |\sqrt{1-\lambda (p)}} .
\end{equation*}
\qed

\begin{notation}\label{not:notation} We will use the following miscellanea of  notations and definitions.
\begin{enumerate}

\item   We will set
\begin{align}
  \label{eq:notation1} \text{$\mathbf{e}(\omega):=\mathbf{E}( \phi _{\omega}) $,  $   \mathbf{q}(\omega):=\mathbf{Q}( \phi _{\omega}) $ and $\mathbf{d}(\omega):=\mathbf{e}(\omega)+\omega \mathbf{q}(\omega)$.}
\end{align}

\item For $z\in \C$ we will use $z_1= \Re z $  and  $z_2= \Im  z $. 

\item Like in the theory of     Kowalczyk et al.   \cite{KMM3},      we   consider constants  $A, B,  \epsilon , \delta  >0$ satisfying
 \begin{align}\label{eq:relABg}
\log(\delta ^{-1})\gg\log(\epsilon ^{-1}) \gg     A  \gg    B^2\gg B  \gg 1.
 \end{align}
 Here we will take $ A\sim B^3$, see Sect. \ref{sec:smooth1} below,  but in fact  $ A\sim B^n$ for  any   $n>2$ would make no difference.

 \item The notation    $o_{\varepsilon}(1)$  means a constant   with a parameter $\varepsilon$ such that
 \begin{align}\label{eq:smallo}
 \text{ $o_{\varepsilon}(1) \xrightarrow {\varepsilon  \to 0^+   }0.$}
 \end{align}
\item For  $\kappa \in (0,1)$    fixed in terms of $p$  and  small enough, we consider
\begin{align}\label{eq:l2w}&
\|  {\eta} \|_{  L ^{p,s}} :=\left \|\< x \> ^s \eta \right \|_{L^p(\R )}   \text{  where $\< x \>  := \sqrt{1+x^2}$,}\\&
\| \eta  \|_{  { \Sigma }_A} :=\left \| \sech \(\frac{2}{A} x\) \eta '\right \|_{L^2(\R )} +A^{-1}\left \|    \sech \(\frac{2}{A} x\)  \eta   \right\|_{L^2(\R )}  \text{ and} \label{eq:normA}\\& \|  {\eta} \|_{ \widetilde{\Sigma} } :=\left \| \sech \( \kappa \omega _0 x\)   {\eta}\right \|_{L^2(\R )} . \label{eq:normk}
\end{align}
\item We set
\begin{align}\label{eq:C+-}
 \C _\pm := \{ z\in \C : \pm \Im z >0   \} .
\end{align}

\item The point $\im \omega$ is a resonance for $\mathcal{L}_\omega$ if there exists a nonzero $ v\in L^\infty (\R , \C ^2)$
such that $\mathcal{L}_\omega v= \im \omega v$. Notice that for $p=3$ the point $\im \omega$ is a resonance.

\item Given two Banach spaces $X$ and $Y$ we denote by $\mathcal{L}(X,Y)$ the space of continuous linear operators from $X$ to $Y$. We write $\mathcal{L}(X ):=\mathcal{L}(X,X)$.

\item We have denoted by $P_c$ the projection on the space \eqref{eq:dirsum11} associated to the spectral decomposition
\eqref{eq:dirsum1}  of the operator $\mathcal{L}_{\omega}$. Later in \eqref{eq:opH} we will introduce an operator $H_{\omega}$ which is an equivalent to $\mathcal{L}_{\omega}$ and  obtained from $\mathcal{L}_{\omega}$ by a simple conjugation. By an abuse of notation we will continue to denote by $P_c$ the analogous spectral projection to the continuous spectrum component, only of $H_\omega$ this time.

\item We have the following elementary formulas,
\begin{align} \label{eq:derf1}&
  Df(u)X=\left . \frac{d}{dt} \(  |u+tX | ^{p-1}\( u+tX \) \)  \right | _{t=0}= |u  | ^{p-1}X + (p-1) |u  | ^{p-3}u \< u, X\> _\C \text{ and}\\&  D^2f(u)X^2= \left . \frac{d}{dt}    Df(u+tX )X  \right | _{t=0} \nonumber\\& = 2(p-1) |u  | ^{p-3} X \< u, X\> _\C + (p-1) |u  | ^{p-3}u |X|^2  + (p-1) (p-3)  |u  | ^{p-5}u\< u, X\> _\C ^2. \label{eq:derf2}
\end{align}

\item Following  the framework in Kowalczyk et al. \cite{KMM3} we   fix an even function $\chi\in C_c^\infty(\R , [0,1])$ satisfying
\begin{align}  \label{eq:chi} \text{$1_{[-1,1]}\leq \chi \leq 1_{[-2,2]}$ and $x\chi'(x)\leq 0$ and set $\chi_C:=\chi(\cdot/C)$  for a  $C>0$}.
\end{align}

\item  We denote by $H^ 1_\rad (\R )= H^ 1_\rad (\R ,\C ) $ the space  of even functions in $H^ 1  (\R ,\C )$.

 \end{enumerate}

\end{notation} \qed

  The group  $e^{t \mathcal{L}_{\omega }}$ is well defined in $L^2  (\R , \C ^2)$, leaves invariant $L^2 (\R , \R ^2)$ and the terms of the direct sums in \eqref{eq:dirsum} and \eqref{eq:dirsum1}.
The following result is an immediate consequence of a Proposition 8.1 in Krieger and Schlag \cite{KrSch} if $p\in \widetilde{\mathbf{F}}$ in \eqref{eq:deftildeF},
since then $\mathcal{L}_{\omega } $ is    admissible in the sense indicated in \cite{KrSch}.
\begin{proposition}\label{prop:KrSch} For any fixed $s>3/2$   there is a  constant $C_\omega $
such that
\begin{align} \label{eq:KrSch}
  \| P_c e^{t \mathcal{L}_{\omega }} : L ^{2,s}  (\R , \C ^2) \to  L ^{2,-s} (\R , \C ^2)  \| \le C_\omega \< t \> ^{-\frac{3}{2}} \text{  for all $t\in \R$.}
\end{align}

\end{proposition}\qed

 For the proof of the following, we refer to  \cite{CM243}.
\begin{proposition} \label{lem:smoothest1} Let   $p\in \widetilde{\mathbf{F}}$.  For  $s>3/2$ and $\tau >1/2$ there exists a constant $C>0 $ such that
 \begin{align}&   \label{eq:smoothest11}   \left \|   \int   _{0} ^{t   }e^{  (t-t')  \mathcal{L}_{\omega }}P_c(\omega )g(t') dt' \right \| _{L^2( \R ,L^{2,-s} (\R ))  } \le C  \|  g \| _{L^2( \R , L^{2,\tau}  (\R ) ) } \text{ for all $g\in  L^2( \R , L^{2,\tau} (\R ) )$}.
\end{align}
\end{proposition}

For the proof of the following, we refer to  \cite{CM24D1}.

\begin{proposition}[Kato smoothing]\label{lem:smooth111}  Let   $p\in \widetilde{\mathbf{F}}$. For any $\omega$ and for any $s>1$  there exists $c>$   such that
\begin{equation}\label{eq:smooth111}
   \|   e^{\im t \mathcal{L}_\omega}P_c u_0 \| _{L^2(\R , L ^{2,-s}(\R   ))} \le c \|  u_0 \| _{L^2(\R  )}.
\end{equation}
\end{proposition}

\section{Refined profile, modulation, continuation argument and proof of Theorem \ref{thm:asstab}}\label{sec:mod}

It is well known,  see Weinstein \cite{W2},  that
\begin{align}
  \label{eq:mangs} \mathcal{S}=\left \{ e^{\frac \im 2 v  x
+\im \vartheta  } \phi _\omega
(\cdot  - D)    : \vartheta , D, v \in \R , \omega >0     \right \}
\end{align}
is a symplectic  submanifold  of  $L ^2  (\R, \C ) $.
We set
 \begin{align}
   \label{eq:refprof} \phi[  {z}] = \phi   + \widetilde{\phi}[  {z}]  \text{ with } \widetilde{\phi}[  {z}] :=  (z \xi  +\overline{z} \overline{\xi}   )     .
 \end{align}
For
functions
\begin{align}
  \label{eq:tildetheta} \widetilde{\Theta}_\mathcal{R}[ z] =(\widetilde{\vartheta}_\mathcal{R}[ z], \widetilde{\omega}_\mathcal{R} [ z], 0 ,0 , \widetilde{z}_\mathcal{R} [ z] )
\end{align}
 in part  to be determined below we introduce
\begin{align*}&
  \widetilde{ z}[ z]=\widetilde{z}_0 [ z]+ \widetilde{z}_\mathcal{R}[ z] \text{ with } \widetilde{z}_0[  z]=\im   \lambda (p)   z\\&\widetilde{\vartheta}[ z]=1 +  \widetilde{\vartheta}_\mathcal{R}[ z] \text{  and }\widetilde{ \omega}[ z]= \widetilde{\omega}_\mathcal{R}[ z] , \quad \widetilde{ v}[ z]= 0 \text{ and } \widetilde{ D}[ z]= 0.
\end{align*}

\begin{proposition}\label{prop:refpropf}
     There exist constants $C_1,\delta _1 >0$ and  a map $\widetilde{\Theta}_\mathcal{R} [  \cdot ]  : D _{\C } (0, \delta _1) \to  \R ^4\times \C  $   like in   \eqref{eq:tildetheta}    s.t. the following   holds:

\begin{description}
  \item[i]  we have
  \begin{align}\label{eq:estpar}
 \left |\widetilde{\Theta}_\mathcal{R} [  {z}] \right |  \lesssim |z|^2;
\end{align}
  \item[ii]   if we set
\begin{align}\label{eq:phi_pre_gali}
R[  {z}]:= \partial ^2_x\phi  [  {z}]+ f(\phi [  {z}]) - \widetilde{\vartheta}\phi [  {z}]+ \im \widetilde{\omega}\Lambda_p\phi [  {z}] +\im  D_{z}\phi  [  {z}]\widetilde{z }   ,
\end{align}
we have for some fixed $ \kappa >0$
\begin{align}&
 \| \cosh \( \kappa  x  \)  {R}   [  {z}]\| _{L^2(\R ) }\lesssim |z|^2 ,\label{estR}
\end{align}
with furthermore the following orthogonality conditions, for $z_1=\Re z$ and $z_2=\Im z$,
\begin{align}\label{R:orth} &
\< \im {R}[  {z}], \phi [  {z}]\>=\< \im {R}[  {z}],\im \Lambda_p \phi [  {z}]\>  = \< \im {R}[  {z}], x\phi [  {z}]\>   =\< \im {R}[  {z}],\im \partial _x\phi [  {z}]\>  \\&  =\< \im {R}[  {z}],\im \partial_{z_{ j}}\phi[  {z}]\> = 0,\text{   for all $j=1,2 $;} \nonumber
\end{align}
  \item[iii]  we have  $\widetilde{\Theta}_\mathcal{R} [  \cdot ]  \in C^1(D _{\C } (0, \delta _1) )$;  \item[iv] $\widetilde{\Theta}_\mathcal{R} [  \cdot ]$ is twice differentiable in $z=0$.
\end{description}

 \end{proposition}
\proof  From  \eqref{eq:static}   and
\begin{align*}
     D_{z}\widetilde{\phi}  [  {z}]\widetilde{z } _0 =    \mathcal{L}  \widetilde{\phi}  [  {z}] =
 -\im \( -\partial _x ^2 \widetilde{\phi}  [  {z}]  +  \widetilde{\phi}  [  {z}] - Df(\phi   )\widetilde{\phi}  [  {z}] \)
\end{align*}
we obtain
\begin{align}\label{eq:defhatR1}& \widehat{R}   [  {z}] = \partial _x ^2 \widetilde{\phi}  [  {z}]   +  f(\phi [  {z}]) - {\phi}  [  {z}] -
 \im  D_{z} {\phi}  [  {z}]\widetilde{z } _0    \text{ where }
    \\& \label{eq:defhatR2}\widehat{R}   [  {z}] := f( \phi +   \widetilde{ \phi} [  {z}]) -  f(\phi )          - Df(\phi   )\widetilde{\phi}  [  {z}] .
\end{align}
\begin{claim}\label{claim:esthatR} We have
\begin{align} \label{eq:esthatR}
   |\widehat{R}   [  {z}]|  \lesssim |z| ^2  e^{|x| \( 2-p-2\sqrt{1-\lambda (p)} \)  } .
\end{align}
\end{claim}
\proof By Lemma \ref{lem:deceig} have $| \widetilde{\phi}  [  {z}]|\lesssim |z| e^{-|x|\sqrt{1-\lambda (p)}} $. We distinguish two case.

\begin{case}\label{case:one}
Let $| \widetilde{\phi}  [  {z}]| \leq  |\phi|/2$.   By $ |\phi| \sim e^{-|x|}$ and $  \left | f''\( \phi +   \tau \widetilde{ \phi} [  {z}]  \) \right |  \sim e^{-(p-2)|x|} $, for $ \tau\in (0, 1) $ we have
\begin{align*}
  |\widehat{R}   [  {z}]|  \lesssim  \sup_{\tau\in(0, 1) }\left | f''\( \phi +   \tau \widetilde{ \phi} [  {z}]  \) \right |  |z| ^2 e^{-2|x|\sqrt{1-\lambda (p)}}     \lesssim     |z| ^2  e^{|x| \( 2-p-2\sqrt{1-\lambda (p)} \)  }
\end{align*}
\end{case}

\begin{case}\label{case:two}
Let $| \widetilde{\phi}  [  {z}]| >  |\phi|/2$.
Then
\begin{align}\label{eq:case:two1}
     e^{-|x|} \lesssim  | \widetilde{\phi}  [  {z}]| \lesssim |z|  e^{-|x| \sqrt{1-\lambda (p)}}   \Longrightarrow  e^{ |x| \(-1 +\sqrt{1-\lambda (p)} \)  } \lesssim |z|
\end{align}
and
\begin{align}\label{eq:case:two2}
   |\widehat{R}   [  {z}]|    & \lesssim  |\widetilde{ \phi} [  {z}]| ^p \lesssim    |z| ^p   e^{ -p|x| \sqrt{1-\lambda (p)}    }   \lesssim   |z| ^p   e^{ -p|x| \sqrt{1-\lambda (p)}    }    |z|  ^{2-p}  e^{ |x| \((2-p) -(2-p)\sqrt{1-\lambda (p)} \)  } \\& = |z| ^2   e^{|x| \( 2-p-2\sqrt{1-\lambda (p)} \)  }  . \nonumber
\end{align}
\end{case}
\qed

 Now we set
\begin{align*}
R[  {z}]=\widehat{R}   [  {z}]   - \widetilde{\vartheta}_\mathcal{R} \phi [  {z}]+ \im \widetilde{\omega}_\mathcal{R} \Lambda_p\phi [  {z}]   -\im \widetilde{D}_\mathcal{R} \partial _x \phi  [  {z}]-\frac{\widetilde{v}_\mathcal{R} x}{2} \phi  [  {z}]+\im  D_{z}\phi  [  {z}]\widetilde{z }_\mathcal{R},
\end{align*}
and impose that the     orthogonality conditions  \eqref{R:orth} be satisfied.
 This  yields
\begin{align}\label{eqsysta}
   \mathbf{A} [z]  \widetilde{\Theta } _\mathcal{R}^\intercal  =   \begin{pmatrix}
\<  \widehat{{R}}[z], \im \phi [  {z}]\> \\   \<  \widehat{{R}}[z], \Lambda_p\phi [  {z}]\>  \\   \<  \widehat{{R}}[z], \partial _x\phi [  {z}]\>   \\   \<  \widehat{{R}}[z], \im x\phi [  {z}]\>  \\   \< \widehat{{R}}[z] , D _ z  \phi\>     \end{pmatrix}    \text{ for an unknown } \widetilde{\Theta } _\mathcal{R}=
\(
 \widetilde{\vartheta } _\mathcal{R} ,  \widetilde{\omega}_\mathcal{R} , \widetilde{D}_\mathcal{R} , \widetilde{v}_\mathcal{R} ,  \widetilde{z } _\mathcal{R}\) .
\end{align}
where \small
 \begin{align}\label{eq:matrixa}
 & \mathbf{A}[z] =\\&  \begin{pmatrix}
 0 &  -\< \Lambda _p \phi [  {z}],   \phi [  {z}]\>  &   0 &  0  &   -\< D _ z  \phi [  {z}],   \phi [  {z}]\>   \\  \< \phi [  {z}], \Lambda _p \phi [  {z}] \> &  0  &   \cancel{ \< \im \partial _x \phi [  {z}],   \Lambda _p \phi[  {z}] \> } &  \cancel {\<  \frac{x}{2} \phi [  {z}],   \Lambda _p \phi  [  {z}] \>  } &  -\< \im  D _ z  \phi [  {z}],  \Lambda _p \phi [  {z}] \>   \\  0   &
 \cancel{-\<\im \Lambda _p \phi [  {z}],   \partial _x\phi [  {z}]\>}  &    0  & -\frac{1}{4}\|  \phi [  {z}]\| _{L^2}^2  &   \cancel{ -\< \im  D _ z  \phi [  {z}],  \partial _x\phi [  {z}] \> } \\  0 &  \cancel{-\<  \Lambda _p \phi [  {z}],     x\phi  [  {z}]\> } &   -\frac{1}{2}\|  \phi [  {z}]\| _{L^2}^2 & 0 &    \cancel{-\<    D _ z  \phi [  {z}],   x\phi   [  {z}]\>}\\ \< \phi [  {z}],   D _ z  \phi  [  {z}]\>   &  -\<\im  \Lambda _p \phi [  {z}] ,    D _ z  \phi [  {z}]\>  &   \cancel{ \< \im \partial _x \phi [  {z}],    D _ z  \phi  [  {z}]\> } & \cancel{ \<  \frac{x}{2} \phi [  {z}],   D _ z  \phi [  {z}]\>   } &    -\<  \im  D _ z  \phi [  {z}],   D _ z  \phi  [  {z}] \>
\end{pmatrix}     \nonumber
\end{align}
\normalsize
where the cancelled terms are null because    by \eqref{eq:sol} and Lemma \ref{lem:simandmult} the function $ \phi [  {z}]$ is even in $x$ and so
\begin{align*}
  \< \im \partial _x \phi [  {z}],   \Lambda _p \phi[  {z}] \> &=  \< \im \partial _x \phi [  {z}],   \Lambda _p \phi[  {z}] \>  =        \< \im  D _ z  \phi [  {z}]  \Xi ,  \partial _x\phi [  {z}] \> \\& =  \< \im  D _ z  \phi [  {z}]  \Xi , \im \partial _x \phi [  {z}] \> = 0 \text{ for any  $\Xi \in \C$.}
\end{align*}
In particular, near $z=0$   we have
 \begin{align} \label{eq:matrixa1}
 & \mathbf{A}[z] = \begin{pmatrix}
 0 &     \frac{p-5}{p-1} q(1)  &   0 &  0  &  O(z)   \\  \frac{5-p}{p-1} q(1) &  0  &  0  &  0 &  O(z)  \\  0   &0 &    0  & -\frac{1}{4}\|  \phi [  {z}]\| _{L^2}^2  &   0  \\  0 &  0  &   -\frac{1}{2}\|  \phi [  {z}]\| _{L^2}^2 & 0 &   0
 \\O(z)  &  O(z)  &  0  & 0 &    J
\end{pmatrix}   + O(z^2) .
\end{align}
So we get \small
\begin{align*}&
   \widetilde{\Theta} _{\mathcal{R}}[z] =\mathbf{A}^{-1}[z]  v[z] \text{  with} \\& v[z]:= \(
\<  \widehat{{R}}[z], \im \phi [  {z}]\> , \<  \widehat{{R}}[z], \Lambda_p\phi [  {z}]\>  , \cancel {\<  \widehat{{R}}[z], \partial _x\phi [  {z}]\> }   ,  \cancel { \<  \widehat{{R}}[z], \im x\phi [  {z}]\>  },  \< \widehat{{R}} [z], D _ z  \phi\>    \) ^\intercal
\end{align*}
\normalsize
where the two cancelled terms are null because $\widehat{{R}}[z]$ is an even function in $x$.    This immediately implies from  system \eqref{eqsysta} that
\begin{align*}
  \widetilde{D}_\mathcal{R} = \widetilde{v}_\mathcal{R}\equiv 0 .
\end{align*}
Properties \textbf{i}, \textbf{iii } and \textbf{iv} need to be checked for $ v[z]$.
Then
\begin{align*}
  |v[z]|\lesssim |z|^2
\end{align*}
follows immediately from  \eqref{eq:esthatR}.  By dominated convergence,   $  v[\cdot ] \in C^0 \( D _{\C } (0, \delta _1)  \). $
This implies \eqref{eq:estpar}  and   $  \widetilde{\Theta}_{\mathcal{R}}[\cdot ] \in C^0 \( D _{\C } (0, \delta _1)  \). $

\noindent Next we consider (notice that $z\to \widetilde{\phi} [z] $  is linear)
\begin{align*}
  D_z \widehat{R} = \(  f'(\phi + \widetilde{\phi} [z]) -  f' (\phi  ) \) D_z \widetilde{\phi} [z]  =   \(  f'(\phi + \widetilde{\phi} [z]) -  f' (\phi  ) \) D_z \widetilde{\phi} [0] .
\end{align*}
 \begin{claim}\label{claim:derbound}
    We have the pointwise estimate
    \begin{align} \label{eq:derbound1}
      |\(  f'(\phi + \widetilde{\phi} [z]) -  f' (\phi  ) \) D_z \widetilde{\phi} [0]| \lesssim   |z|   e^{|x| \( 2-p-2\sqrt{1-\lambda (p)} \)  } .
    \end{align}
 \end{claim}
 \proof  For  $| \widetilde{\phi}  [  {z}]| \leq  |\phi|/2$,  by $ |\phi| \sim e^{-|x|}$ and $  \left | f''\( \phi +   \tau \widetilde{ \phi} [  {z}]  \) \right |  \sim e^{-(p-2)|x|} $, for $ \tau\in (0, 1) $ we have \small
\begin{align*}
   |\(  f'(\phi + \widetilde{\phi} [z]) -  f' (\phi  ) \) D_z \widetilde{\phi} [0]|  \lesssim  \sup_{\tau\in(0, 1) }\left | f''\( \phi +   \tau \widetilde{ \phi} [  {z}]  \) \right |  |z|   e^{-2|x|\sqrt{1-\lambda (p)}}     \lesssim     |z|   e^{|x| \( 2-p-2\sqrt{1-\lambda (p)} \)  } .
\end{align*}
 \normalsize
 For  $| \widetilde{\phi}  [  {z}]|  \gtrsim  |\phi| $  like in  \eqref{eq:case:two1}--\eqref{eq:case:two2}
\begin{align*} &
  |\(  f'(\phi + \widetilde{\phi} [z]) -  f' (\phi  ) \) D_z \widetilde{\phi} [0]|   \lesssim    |     \widetilde{\phi} [z]   | ^{p-1}  e^{- |x|\sqrt{1-\lambda (p)}}    \lesssim |z|^{p-1}
   e^{- p|x|\sqrt{1-\lambda (p)}} \\& \lesssim  |z|^{p-1}
   e^{- p|x|\sqrt{1-\lambda (p)}} |z|   ^{2-p}  e^{ |x| \((2-p) -(2-p)\sqrt{1-\lambda (p)} \)  }      \lesssim  |z|    e^{|x| \( 2-p-2\sqrt{1-\lambda (p)} \)  }.\nonumber
\end{align*}
 \qed

 Claim \ref{claim:derbound} allows to conclude that we can use differentiation under the integral and that  $  v[\cdot ] \in C^1 \( D _{\C } (0, \delta _1)  \)  $  with
 \begin{align*}
  |D_zv[z]|\lesssim |z|
\end{align*}
 This in turn implies   $  \widetilde{\Theta}_{\mathcal{R}}[\cdot ] \in C^1 \( D _{\C } (0, \delta _1)  \)  $ with a similar bound.

 \noindent We consider now, formally
 \begin{align*}
  D_z ^2\widehat{R} [z]=    f''(\phi + \widetilde{\phi} [z])   ( D_z \widetilde{\phi} [z] )^2 =  f''(\phi + \widetilde{\phi} [z])   ( D_z \widetilde{\phi} [0] )^2.
\end{align*}
Now we write
\begin{align*}
    D_z  \widehat{R} [z] -   D_z ^2\widehat{R} [0]z = \(    f'(\phi + \widetilde{\phi} [z]) -  f' (\phi  )   -   f''(\phi  )    \widetilde{\phi} [z] \)     D_z \widetilde{\phi} [0]    =:I.
\end{align*}
For  $| \widetilde{\phi}  [  {z}]| \leq  |\phi|/2$, by $|f'''(\phi +\tau s \widetilde{\phi} [z])| \lesssim |f'''(\phi  )|$  for $s,\tau \in [0,1]$ and $\varepsilon _0>0$ small,
\begin{align*}
  |I|&\lesssim  \int _{[0,1]^2} \tau |f'''(\phi +\tau s \widetilde{\phi} [z]) | ds d\tau     |\widetilde{\phi} [z]| ^{1+\varepsilon _0}     \frac{|\widetilde{\phi} [z]|^{1-\varepsilon _0}}{ \phi ^{1-\varepsilon _0}}   \phi ^{1-\varepsilon _0}  |D_z \widetilde{\phi} [0] | \\&  \lesssim   |z| ^{1+\varepsilon _0}        e^{|x|  \( 2+ \varepsilon _0 -p      -2 \sqrt{1-\lambda (p)}   \)} .
\end{align*}
For   $| \widetilde{\phi}  [  {z}]|  \gtrsim  |\phi| $ we similarly have
\begin{align*}
  |I| &\lesssim    | \widetilde{\phi}  [  {z}]| ^{p-1} | D_z \widetilde{\phi} [0]| \lesssim  |z|  ^{p-1} e^{-p |x| \sqrt{1-\lambda (p)}} \\& \lesssim   |z|  ^{p-1} e^{-p |x| \sqrt{1-\lambda (p)}}
  |z|^{2-p+\varepsilon _0} e ^{|x| (2-p+\varepsilon _0) \(   1 - \sqrt{1-\lambda (p)}\)     } \le  |z|  ^{1+\varepsilon _0}  e^{|x|  \( 2+ \varepsilon _0 -p      -2 \sqrt{1-\lambda (p)}   \)} .
\end{align*}
 Now let us consider the following component of $v[z]$ (the other can be treated similarly)
 \begin{align*}
   & \<  \widehat{{R}}[z], \im \phi [  {z}]\> = \<  \widehat{{R}}[z], \im \phi  \>  +\<  \widehat{{R}}[z], \im \widetilde{ \phi} [  {z}]   \> =II_1+II_2.
 \end{align*}
Then
\begin{align*}
   \lim _{z\to 0}  \frac{\<   D_z  \widehat{R} [z] , \im \phi  \>   - \<   D_z ^2\widehat{R} [0]z  , \im \phi  \>    }{|z|} =0
\end{align*}
which implies the differentiability of $II_1$ at $z=0$. The same is true for $II_2$, with differential 0 at $z=0$.  Since the other components of $v[z]$  can be treated similarly, this  gives the 2nd order differentiability of $v[z]$, and thus also of  $\widetilde{\Theta}_\mathcal{R} [  z ]$,    at $z=0$.

\qed

If $(\mathcal{D}_{\sqrt{\omega}} u)(x):=  \omega ^{\frac 1{p-1}}  u(\sqrt{\omega} x) $ we get the following.
\begin{lemma}
  \label{lem:dil} Set $\phi[\omega,z](x) := \omega ^{\frac 1{p-1}} \phi[z](\sqrt{\omega}x)$. Then
  \small\begin{align} \label{eq:refomega} &
R[ \omega ,  {z}] = \partial ^2_x\phi  [  \omega ,{z}]+ f(\phi [ \omega , {z}]) - \widetilde{\vartheta}[ \omega , {z}]\phi [ \omega , {z}]+ \im \widetilde{\omega}[ \omega , {z}] \partial _\omega \phi [  \omega ,{z}]   +\im  D_{z}\phi  [ \omega , {z}]\widetilde{z }[ \omega , {z}]
\end{align}\normalsize
for  \begin{align} &
 \label{def:phiomegaz}
R[\omega,z](x):=\omega^{\frac{p}{p-1}}R[z](\sqrt{\omega}x), \\ &
\widetilde{\vartheta}[\omega,z]:=\omega \widetilde{\vartheta}[z], \quad \widetilde{\omega}[\omega,z]=  \widetilde{\omega}[z] \text{ and }
\widetilde{z}[\omega,z]= \omega \widetilde{z}[z].\nonumber
\end{align}
Furthermore, we have the following orthogonality  properties,
\begin{align} \nonumber &
\< \im R[ \omega  ,  {z}], \phi [ \omega  ,  {z}] \>=\< \im R[ \omega  ,  {z}] ,\im \partial _\omega  \phi [ \omega  ,  {z}] \>  \\&   = \< \im R[ \omega  ,  {z}] , x\phi [ \omega  ,  {z}] \>   =\< \im R[ \omega , {z}] ,\im \partial _x\phi [ \omega  ,  {z}] \> \label{R:orthomega}  \\&  =\< \im R[ \omega  ,  {z}] ,\im \partial_{z_{ j}}\phi[ \omega  ,  {z}]\> = 0,\text{   for all $j=1,2 $.} \nonumber
\end{align}
\end{lemma}
\proof Notice that
\begin{align*}
  \mathcal{D}_{\sqrt{\omega}}  \partial _x =  \omega ^{-1/2} \partial _x  \mathcal{D}_{\sqrt{\omega}} \quad , \quad   \mathcal{D}_{\sqrt{\omega}}   x =  \omega ^{ 1/2}  x  \mathcal{D}_{\sqrt{\omega}} \text{ and }
  \mathcal{D}_{\sqrt{\omega}} f(\phi [   {z}]) = \omega ^{   - 1} f(\phi [ \omega , {z}]).
\end{align*}
So that applying  $\mathcal{D}_{\sqrt{\omega}}$ to equation \eqref{eq:phi_pre_gali}   and using  $[\mathcal{D}_{\sqrt{\omega}}, \Lambda _p]=0$
we obtain \small
\begin{align*}
 \mathcal{D}_{\sqrt{\omega}} R[  {z}] &=   \omega ^{-1 }  \partial ^2_x\phi  [ \omega , {z}]+ \omega ^{   -1} f(\phi [ \omega , {z}])  - \widetilde{\vartheta}[z] \phi [ \omega , {z}]+ \im \widetilde{\omega}\Lambda_p\phi [ \omega ,  {z}]   +\im  D_{z}\phi  [ \omega ,  {z}]\widetilde{z }  .
\end{align*}\normalsize
Multiplying  by $\omega $, using $\partial_{\omega}\phi [ \omega , {z}] =  \omega ^{-1} \Lambda_p  \phi [ \omega , {z}] $ and  the definitions in \eqref{def:phiomegaz} we obtain  \eqref{eq:refomega}.

The  orthogonality  properties  \eqref{R:orthomega} follow immediately from \eqref{R:orth}.
\qed

 \begin{lemma}  \label{lem:boost} Set $\phi[\omega , v,z]  :=  e^{\frac{\im }{2} vx}\phi[\omega  ,z] $. Then \begin{align} \label{eq:refboost} &R[ \omega , v,  {z}] = \partial ^2_x\phi  [  \omega , v,{z}]+ f(\phi [ \omega ,v,  {z}]) - \widetilde{\vartheta}[ \omega ,v,  {z}]\phi [ \omega  ,v , {z}]+ \im \widetilde{\omega}[ \omega ,v,  {z}] \partial _\omega \phi [  \omega ,v, {z}]\\&  +\im  D_{z}\phi  [ \omega , v, {z}]\widetilde{z }[ \omega ,v,  {z}] -\im v \partial _x \phi  [ \omega ,v,  {z}] .\nonumber\end{align}for  \begin{align} & \label{def:phiboost}R[\omega, v,z] :=e^{\frac{\im }{2} vx}  R[ \omega ,  {z}] ,     \quad \widetilde{\omega}[\omega,v, z]=\widetilde{\omega}_{\mathcal{R}} [\omega,v, z]  \text{ with }    \widetilde{\omega}_{\mathcal{R}} [\omega,v, z] =\widetilde{\omega}_{\mathcal{R}}[  z]  , \\   
 &\widetilde{z}[\omega , v,z]= \widetilde{z}_0[\omega , v,z]     +   \widetilde{z}_{\mathcal{R}}[\omega , v,z]        \text{ with } \widetilde{z}_0[\omega , v,z]=\im \omega \lambda (p) z  \text{ and }  \widetilde{z}_{\mathcal{R}}[\omega , v,z] =   \omega \widetilde{z}_{\mathcal{R}}[  z]  ,  \nonumber  \\& \widetilde{\vartheta }[ \omega , v , {z}] =\omega -\frac{v^2}{4}  + \widetilde{\vartheta} _{\mathcal{R}}[ \omega  ,v, {z}] \text{ with }   \widetilde{\vartheta} _{\mathcal{R}}[ \omega  ,v, {z}] =  \omega  \widetilde{\vartheta} _{\mathcal{R}}[   {z}] .\nonumber\end{align}Furthermore, we have the following orthogonality  properties,\begin{align} \nonumber &\< \im R[ \omega , v,  {z}], \phi [ \omega , v,  {z}] \>=\< \im R[ \omega , v,  {z}] ,\im \partial _\omega  \phi [ \omega , v,  {z}] \>  \\&   = \< \im R[ \omega , v,  {z}] , x\phi [ \omega , v,  {z}] \>   =\< \im R[ \omega , v,  {z}] ,\im \partial _x\phi [ \omega , v,  {z}] \> \label{R:orthboost}  \\&  =\< \im R[ \omega , v,  {z}] ,\im \partial_{z_{ j}}\phi[ \omega , v,  {z}]\> = 0,\text{   for all $j=1,2 $.} \nonumber\end{align}\end{lemma}\proof  We multiply \eqref{eq:refomega}  by $e^{\frac{\im }{2} vx}$ and use\begin{align*}&    e^{\frac{\im }{2} vx} \partial ^2_x\phi  [  \omega ,{z}] = \partial ^2_x\phi  [  \omega , v,{z}] +\frac{v^2}{4}\phi  [  \omega , v,{z}] -\im v \partial _x \phi  [  \omega , v,{z}]\end{align*}obtaining the following, which yields  \eqref{eq:refboost},\begin{align*} & e^{\frac{\im }{2} vx}  R[ \omega ,  {z}] =  \partial ^2_x\phi  [  \omega , v,{z}] +\frac{v^2}{4}\phi  [  \omega , v,{z}] -\im v \partial _x \phi  [  \omega , v,{z}]  + f(\phi [ \omega , v , {z}])\\&  - \widetilde{\vartheta}[ \omega , {z}]\phi [ \omega , v , {z}]+ \im \widetilde{\omega}[ \omega , {z}] \partial _\omega \phi [  \omega , v ,{z}] +\im  D_{z}\phi  [ \omega , {z}]\widetilde{z }[ \omega , v , {z}]  .\end{align*}The  orthogonality  properties  \eqref{R:orthboost} follow immediately from \eqref{R:orth}.\qed

Let now
\begin{align*}&    \widetilde{v}[\omega,v, z]=0    \text{ and }  \widetilde{D}[ \omega , v , {z}] =  v
\end{align*}
and
\begin{align}&
\Theta = \( \vartheta , \omega , D, v,  z\) , \nonumber \\
&
 \label{eq:tildethetab}
   \widetilde{\Theta}_\mathcal{R}[ \Theta ]  =(\widetilde{\vartheta}_\mathcal{R}[  \omega , v,{z} ], \widetilde{\omega}_\mathcal{R} [  \omega , v,{z} ],0 , 0 , \widetilde{z}_\mathcal{R} [  \omega , v,{z} ] ) \text{ and }\\&  \widetilde{\Theta} [ \Theta ] =(\widetilde{\vartheta} [  \omega , v,{z} ], \widetilde{\omega}  [  \omega , v,{z} ],0 ,v , \widetilde{z}  [  \omega , v,{z} ] ) \text{ and}\nonumber \\& \phi  [ \Theta ] = e^{ \im \vartheta-\im D  \partial _x }    e ^{\frac{\im  } 2 vx}\phi [ \omega   , {z}].\nonumber
\end{align}
The proof of the following   modulation is standard, see Stuart \cite{stuart}.

\begin{lemma}[Modulation]\label{lem:mod1}  Let $\omega _0 >0$. There exists an $\delta _0 >0$ and
$\Theta     \in C^1(   \mathcal{U} (\omega _0,\delta _0   ), \R   ^4 \times \C ) $
 such that for  any $u \in \mathcal{U} (\omega _0,\delta _0   )$
\small \begin{align}\label{61} &   \eta (u):=e^{-\im \vartheta( u) +  D(u) \partial _x } u -  \phi [\Theta (u)]    \text{ satisfies }   \\& \nonumber\< \im   \eta (u), e^{-\im \vartheta( u) +  D(u) \partial _x } \left . D_{\Theta} \phi  [ \Theta ] \right | _{\Theta =\Theta (u)} \Xi   \> = 0,\text{   for all $\Xi \in  \R   ^4 \times \C$.}\end{align} \normalsize
 \normalsize
 Furthermore we have the identities
\begin{align*}& \omega  ( \phi  [ \Theta ])  =\omega  \, , \quad v (  e^{ \im \vartheta _0-  D _0 \partial _x }  u)    = v  (    u)  \, , \quad z(  e^{ \im \vartheta _0-  D _0 \partial _x }  u)    = z (    u)  \, \\&  \vartheta  (  e^{ \im \vartheta _0-  D _0 \partial _x }  u)    = \vartheta  (    u)   + \vartheta _0
  \text{ and} \\&   D  (  e^{ \im \vartheta _0- D _0 \partial _x }  u)    = D  (    u)    + D_0 .
\end{align*}

\end{lemma}
\qed

We  have now the ansatz
\begin{align}\label{eq:ansatz}
  u = e^{\im \vartheta - D   \partial _x}      \(  \phi [\omega ,v  , z ]+ \eta \) =:e^{\im \vartheta - D   \partial _x}     e^{\frac{\im}{2} vx}\(  \phi _\omega + r \)   \text{ where } r= \widetilde{\phi} [\omega ,  z ]+  e^{-\frac{\im}{2} vx}\eta.
\end{align}
We  can prove that if  $u_0\in D _{H^1 (\R )}( \phi _{\omega  _0} , \delta)$ then         there exist   $(\vartheta  , D)=      (\vartheta  (t) , D)(t)     $ such that
\begin{equation}\label{eq:orbsta1}
 \| u- e^{\im \theta -  D   \partial _x } \phi _{\omega  _0}\| _{H^1} \lesssim  \delta \text{  for all values of time.}
\end{equation}
Notice now that
\begin{align}& \label{eq:expen}
 \mathbf{E} \( e^{\frac{\im}{2} vx} u  \) =  \mathbf{E} \(   u  \)+ v \mathbf{P}(u) + \frac{v^2}{4}\mathbf{Q}(u) \text{  and} \\&  \mathbf{P} \( e^{\frac{\im}{2} vx} u  \) =  \mathbf{P} \(   u  \) +    \frac{v }{2}\mathbf{Q}(u). \label{eq:expmom}
\end{align}
In analogy to  Stuart \cite[ pp. 67--68 ]{stuart} we consider
\begin{align*} &  \mathcal{H } (u):= \mathbf{E}(u ) + \( \omega   +  \frac{v^2(0)}{4}\)  \mathbf{Q}(u)-v(0)   \mathbf{P}(u) .
\end{align*}
Entering the ansatz in \eqref{eq:ansatz} we  get
\begin{align*}&  \mathcal{H}( u ) =
  \mathcal{H}( e^{\frac{\im}{2} vx}\(  \phi _\omega + r \) ) = \mathbf{E}( e^{\frac{\im}{2} vx}\(  \phi _\omega + r \)) + \(\omega   + \frac{v^2(0)}{4} \)  \mathbf{Q}( \phi _\omega + r  )-v(0)  \mathbf{P}(e^{\frac{\im}{2} vx}\(  \phi _\omega + r \) ) \\& = \mathbf{E}( \phi _\omega + r ) + \left ( \omega   + \frac{v^2(0)}{4}   + \frac{v^2-2v(0)
v}{4} \right ) \mathbf{Q}(\phi _\omega + r )+(v-v(0)) \mathbf{P}( \phi _\omega + r ) \\& =
\mathbf{E}( \phi _\omega + r) + \left ( \omega    +
\frac{(v-v(0))^2}{4} \right ) \mathbf{Q}( \phi _\omega + r)+(v-v(0))\mathbf{P}(  \phi _\omega + r ).
\end{align*}
Then expanding the above functionals around $\phi _\omega$    we obtain
\begin{align*}
  & H(u)=\mathbf{d}(\omega  ) + \frac{(v-v(0))^2}{4}   \mathbf{q}(\omega  )+\langle \cancel{\nabla \mathbf{E}  (\phi
_\omega  ) + \omega \nabla \mathbf{Q}  (\phi _\omega  )}+(v-v(0)) \nabla \mathbf{P} ( \phi _\omega ),r \rangle
+\\& + \frac 12 \langle \left [  \nabla ^2 \mathbf{E}  (\phi _\omega  ) + \omega    \nabla ^2  \mathbf{Q}  (\phi
_\omega  ) \right ] r,r \rangle + \frac{(v-v(0))^2}{4} \langle
   \nabla \mathbf{Q}  (\phi _\omega  )+ \frac{ 1}{2} \nabla ^2\mathbf{Q} (\phi _\omega  )  r ,r \rangle +\\& +\frac{1}{2}          (v-v(0))
    \langle  \nabla ^2  \mathbf{P }(
\phi _\omega )r,r \rangle +o(\| r\| _{H^1}^2)
\end{align*}
where the cancelled term is null by \eqref{eq:static}. Since $r=  \widetilde{\phi} [\omega   , z ] +  \eta $, we have
\begin{align*} &
   \<  \nabla \mathbf{P} ( \phi _\omega ),r \>  = - \<  \im \partial _x  \phi _\omega  ,r \>  =   - \cancel{\<  \im \partial _x  \phi _\omega  ,  \widetilde{\phi} [\omega   , z ]      \> }  - \<  \im     \partial _x  \phi _\omega  ,    \eta     \>
   =O(|z| \| \eta \| _{\widetilde{\Sigma}})
\end{align*}
where the cancellation is due to the fact that $\partial _x  \phi _\omega$ is odd and $ \widetilde{\phi} [\omega   , z ]$ is even and the last equality follows from the modulation.
 So
 \begin{align*}&
    H(u)=d(\omega  ) + \frac{(v-v(0))^2}{4}   q(\omega )+ \frac 12 \langle \left [ \nabla ^2 \mathbf{E}  (\phi _\omega  ) + \omega    \nabla ^2  \mathbf{Q}  (\phi
_\omega  )  \right ] r,r \rangle  +o(\| r\| _{H^1}^2 + (v-v(0))^2   )
 \end{align*}
 where we use
 \begin{align*}
   \| r\| _{H^1}^2\sim |z|^2 + \| \eta\| _{H^1}^2
 \end{align*}
 which follows by modulation.  Since  we have  \small
 \begin{align*}
   &\mathcal{H} (u_0) =\mathbf{E}(u_0)+   \(\omega   (t) + \frac{v^2(0)}{4} \)
\mathbf{Q}(u_0)-v(0)\mathbf{P}(u_0) \\& =d(\omega (0)) +   \frac 12 \langle \left [
\nabla ^2 \mathbf{E}  (\phi _{\omega (0)}) +            {\omega (0)}   \nabla ^2  \mathbf{Q} (\phi _{\omega (0)} )
\right ] r(0),r(0)\rangle + (\omega  (t) -{\omega(0)}) q ( {\omega (0)}  ) +o(\| r (0)\| _{H^1}^2   )
 \end{align*} \normalsize
 then by $\mathcal{H} (u (t) ) =\mathcal{H} (u_0) $ and   $\mathbf{d}'(\omega )= \mathbf{q}(\omega)$ we get
 \begin{align*}
    &  \mathbf{d}(\omega  ) -\mathbf{ d}(\omega (0)) - (\omega    -{\omega(0)})\mathbf{d}'  ( {\omega (0)} + \frac{(v-v(0))^2}{4}   \mathbf{q}(\omega  )   + \frac 12 \langle \left [  \nabla ^2 \mathbf{E}  (\phi _\omega  ) + \omega    \nabla ^2  \mathbf{Q}  (\phi
_\omega  ) \right ] r,r \rangle  \\& \lesssim  \| r (0)\| _{H^1}^2  \lesssim \delta ^2 .
 \end{align*}
 Since
 \begin{align*}
   \langle \left [  \nabla ^2 \mathbf{E}  (\phi _\omega  ) + \omega    \nabla ^2  \mathbf{Q}  (\phi
_\omega  ) \right ] r,r \rangle \gtrsim \| r  \| _{H^1}^2
 \end{align*}
 and
 \begin{align*}
   \mathbf{d}(\omega  ) -\mathbf{ d}(\omega (0)) - (\omega    -{\omega(0)})\mathbf{d}'  ( {\omega (0)} =  \frac{\mathbf{d}''(\omega (0) )}{2}(\omega    -{\omega(0)})^2 +o\( (\omega    -{\omega(0)})^2\)
 \end{align*}
 with   $\mathbf{d}''(\omega (0) )>0 $  and by $|\omega (0)-\omega _0|\lesssim \delta$     we conclude that
\begin{align}
  \label{eq:oertstab1} |\omega -\omega _0|+ |v-v (0)|+ |z|+ \| \eta \| _{H^1}\lesssim  {\delta } \text{   for all values of time}.
\end{align}

The proof of Theorem \ref{thm:asstab} is mainly  based on the following continuation argument.

\begin{proposition}\label{prop:continuation}
There exists  a    $\delta _0= \delta _0(\epsilon )   $ s.t.\  if
\begin{align}
  \label{eq:main2} \| \eta \| _{L^2(I, \Sigma _A )} +  \| \eta \| _{L^2(I, \widetilde{\Sigma}  )} + \| \dot \Theta - \widetilde{\Theta} \| _{L^2(I  )} + \| z^2 \| _{L^2(I  )}\le \epsilon
\end{align}
holds  for $I=[0,T]$ for some $T>0$ and for $\delta  \in (0, \delta _0)$
then in fact for $I=[0,T]$    inequality   \eqref{eq:main2} holds   for   $\epsilon$ replaced by $   o _{\varepsilon }(1) \epsilon $.
\end{proposition}
Notice that this implies that in fact the result is true for $I=\R _+$.   We will split the proof of Proposition \ref{prop:continuation} in a number of partial results obtained assuming the hypotheses of Proposition  \ref{prop:continuation}.
\begin{proposition}\label{prop:modpar} We have
  \begin{align}&
  \label{eq:modpar1} \|\dot \vartheta -\widetilde{\vartheta} \|  _{L^1(I  )} + \|\dot \omega -\widetilde{\omega}  \|  _{L^1(I  )} + \|\dot D -v  \|  _{L^1(I  )} +  \|\dot v    \|  _{L^1(I  )}\lesssim  \epsilon ^2 ,\\&  \label{eq:modpar2}   \|\dot z -\widetilde{z} \|  _{L^2(I  )} \lesssim  { \delta}\epsilon  ,\\&  \label{eq:modpar3}   \|\dot z \|  _{L^\infty(I  )} \lesssim   {\delta} .
\end{align}

\end{proposition}

\begin{proposition}[Fermi Golden Rule (FGR) estimate]\label{prop:FGR}
 We have
 \begin{align}\label{eq:FGRint}
  \| z^2\|_{L^2(I)}\lesssim  A^{-1/2} \epsilon .
 \end{align}
 \end{proposition}

\begin{proposition}[Virial Inequality]\label{prop:1virial}
 We have
 \begin{align}\label{eq:sec:1virial1}
   \| \eta \| _{L^2(I, \Sigma _A )} \lesssim  A \delta + \| z^2\|_{L^2(I)} + \| \eta \| _{L^2(I, \widetilde{\Sigma}   )} + \epsilon ^2  .
 \end{align}
 \end{proposition}
Notice that $A \delta \ll  A ^{-1}\epsilon ^2 =o _{B^{-1}} (1)  \epsilon ^2 $  in \eqref{eq:sec:1virial1}.

\begin{proposition}[Smoothing Inequality]\label{prop:smooth11}
 We have
 \begin{align}\label{eq:sec:smooth11}
   \| \eta \| _{L^2(I, \widetilde{\Sigma}   )} \lesssim  o _{B^{-1}} (1)  \epsilon    .
 \end{align}
 \end{proposition}

%

\textit{Proof of Theorem \ref{thm:asstab}.}
It is straightforward that Propositions \ref{prop:modpar}--\ref{prop:smooth11}  imply Proposition \ref{prop:continuation} and thus the fact that we can take $I=\R _+$ in all the above inequalities. This in particular implies \eqref{eq:asstab2}. By $z\in L^4(\R  _+)$ and $\dot z \in L^\infty (\R  _+)$ we have  \eqref{eq:asstab3}.  

\noindent We next focus on the limit \eqref{eq:asstab20}. We first rewrite our equation,
entering the  ansatz \eqref{eq:ansatz} in \eqref{eq:nls1} and using   \eqref{eq:refomega}    we obtain
\begin{align}     \nonumber &
    \dot \eta   +\im  \dot \vartheta  \eta      -  \dot D    \partial _x\eta       + e^{-\im \vartheta +D\partial _x}  D_{\Theta} \phi [\Theta ]  \(\dot  \Theta   - \widetilde{\Theta} \)   \\& =  \im   \partial _x^2    \eta  +  \im  \( f( \phi [\omega ,v  , z ] + \eta  )  -  f( \phi  [\omega ,v, z ]  )  \) +\im R[\omega ,v   , z ]    .   \label{eq:nls3}
\end{align}
Then we can proceed like in \cite{CM24D1} considering
  \begin{equation*}
  \begin{aligned}
   \mathbf{a}(t) &:= 2 ^{-1}\|  e^{- a\< x\>}   \eta (t)  \| _{L^2(\R )} ^2
  \end{aligned}
  \end{equation*}
and  by  obvious cancellations and       orbital stability getting  \small
\begin{align}\label{eq:intparts}
  \dot {\mathbf{a}}  &=  \frac{1}{2} \< \left [   e^{- 2a\< x\>} , \im  \partial _x ^2  \right ] \eta , \eta \> +   \frac{1}{2}   \dot D       \< \left [   e^{- 2a\< x\>} ,   \partial _x   \right ] \eta , \eta \>
   -\< e^{-\im \vartheta +D\partial _x}  D_{\Theta} \phi [\Theta ]  \(\dot  \Theta   - \widetilde{\Theta} \)  , e^{- 2a\< x\>} \eta    \>   \\& + \<
  e^{- a\< x\>}  \(     \im  \( f( \phi [\omega ,v  , z ] + \eta  )  -  f( \phi  [\omega ,v , z ]  )  \) +\im R[\omega ,v, z ] \) , e^{- a\< x\>} \eta    \> =O( \delta  ^2 ) \text{ for all times}.\nonumber
\end{align}\normalsize
   Since we already know from \eqref{eq:asstab2} that $\mathbf{a}\in L^1(\R )$, we conclude that
$ \mathbf{a}(t) \xrightarrow{t\to +\infty} 0$.
Notice that the integration by parts in \eqref{eq:intparts} can be made rigorous   considering that
if $u_0\in H^2(\R )$ by the well known regularity result by Kato, see \cite{CazSemi}, we have $\eta \in C^0\( \R , H^2 (\R )\)$ and the above argument is correct and by a standard  density argument the result  can be extended to $u_0\in H^1(\R )$.

We prove \eqref{eq:asstab3}.
Here, notice that by orbital stability we can  take   $a>0$ such that we have the  following, which will be used below,
\begin{equation}\label{eq:uniformity1}
  e^{-2a\<x\>}\gtrsim    \max \{ |\phi[\omega (t),v(t) ,z(t)]| , |\phi[\omega (t),v(t),z(t) ]| ^{p-1} \}  \text {  for all } t\in \R .
\end{equation}
Since $\mathbf{Q}(\phi_{\omega})=\mathbf{q}(\omega)$ is strictly  monotonic in $\omega$, it suffices to show the  $\mathbf{Q}(\phi_{\omega (t)})$ converges as $t\to\infty$.
From the conservation of $\mathbf{Q}$, the exponential decay of $\phi[\omega,v,z]$, \eqref{eq:asstab20} and \eqref{eq:asstab2}, we have
\begin{align}\label{eq:equipart1}
\lim_{t\to \infty}\(\mathbf{Q}(u_0)-\mathbf{Q}(\phi_{\omega (t)})-\mathbf{Q}(\eta (t))\)=0.
\end{align}
Thus, our task is now to prove $\frac{d}{dt}\mathbf{Q}(\eta) \in L^1$, which is sufficient to show the convergence of $\mathbf{Q}(\eta)$.
Now, from \eqref{eq:nls3}, we have
\begin{align*}
\frac{d}{dt}\mathbf{Q}(\eta)&=\<\eta,\dot{\eta}\>
=-\<\eta,e^{-\im \vartheta +D\partial _x}  D_{\Theta} \phi [\Theta ]  \(\dot  \Theta   - \widetilde{\Theta} \)      \>    \\&\quad+\<\eta,  \im  \( f( \phi [\omega ,v  , z ] + \eta  )  -  f( \phi  [\omega ,v  , z ]  )  \)\>+\<\eta, \im R[\omega ,v  , z ]\>\\&
=I+II+III .
\end{align*}
By the bound of the 1st and the 3rd term of \eqref{eq:main2}, we have $I\in L^1(\R _+)$.
$III\in L^1(\R _+)$ follows from \eqref{estR} and $|z|^2\in L^2(\R _+)$. So the key term is $II$.  We partition  for $s\in (0,1)$ the line   where $x$ lives as
\begin{align*}&
   \Omega_{1,t ,s}=\{x\in \R\ |\  | s \eta(t,x)|\leq 2|\phi[\omega (t),v(t) ,z(t)]|\}  \text{ and }\\&  \Omega_{2,t,s }=\R\setminus \Omega_{1,t,s }=\{x\in \R\ |\ | s\eta(t,x)|> 2|\phi[\omega (t),v(t),z(t)]|\} ,
\end{align*}
Then, we have
\begin{align*}
II(t)&=\sum_{j=1,2}   \int_0^1ds \int_{\Omega_{j,t ,s }} \< \im  D  f( \phi[\omega (t),v(t) ,z(t)]+ s\eta (t))  \eta (t) , \eta (t)        \> _\C     \,dx\\& =:II_1(t)+II_2(t).
\end{align*}
For $II_1$,  by    \eqref{eq:derf2}, we have \small
\begin{align*}
|II_1(t)|&\lesssim   \int_0^1ds\int_{\Omega_{1,t,s }}  \left |   \< \im  D  f( \phi[\omega (t),v(t) ,z(t)]+ s\eta (t))  \eta (t) , \eta (t)        \> _\C        \right | dx   \\& \lesssim \int_0^1ds\int_{\Omega_{1,t ,s}}  |\phi[\omega (t),v(t) ,z(t)]  +s \eta (t) | ^{p-1}|\eta (t)|^2\,dx \lesssim \int_{\R} |\phi[\omega (t),v(t),z(t)]|^{p-1}|\eta (t)|^2\,dx \\&  \lesssim \|  e^{- a\<x\>}    \eta (t)\| _{L^2 (\R )}  ^2 \in L^1(\R _+ ).
\end{align*}\normalsize
Turning to $II_2$,  by exploiting $ \<  D  f( s \eta (t)) \eta (t)           ,\im  \eta (t) \> _\C  \equiv 0  $,  which can be easily checked from   \eqref{eq:derf1},
we write \small
\begin{align*}
   II_2(t) &= - \int _{[0,1]  }     ds \int_{\Omega_{2,t,s }} \<  D  f(\phi[\omega (t),v(t) ,z(t)]+ s\eta (t)) \eta (t) -D  f( s\eta (t)) \eta (t)           ,\im  \eta (t) \> _\C      dx \\&  = - \int _{[0,1] ^2 }  d\tau ds  \int_{\Omega_{2,t,s }} \<  D  ^2 f(\tau \phi[\omega (t),v(t),z(t)]+s \eta (t))  (  \phi[\omega (t),v(t),z(t)] ,\eta (t))             ,\im  \eta (t) \> _\C      dx .
\end{align*}\normalsize
Then by  $|s \eta(t,x)|> 2|\phi[\omega (t),v(t),z(t)]|$,   $p\in (1,2]$, \eqref{eq:uniformity1} and orbital stability we have the following, which yields \eqref{eq:asstab3},
\begin{align*}
   |II_2(t)| &\lesssim    \int _{[0,1] ^2 }  d\tau  ds \int_{\Omega_{2,t,s }}       | s \eta (t)  | ^{p-2}  |  \phi[\omega (t),v(t),z(t)]  |   \ |\eta (t) |^2     dx \\& \lesssim   \int _{[0,1]   }    ds \ s^{p-2}  \int_{\Omega_{2,t,s  }}       |  \eta (t)  | ^{p-2}  |  \phi[\omega (t),v(t),z(t)]  | ^{p-1}  |  \eta (t)  | ^{2-p}  \ |\eta (t) |^2     dx \\&\lesssim    \|  e^{- a\<x\>}    \eta (t)\| _{L^2 (\R )}  ^2 \in L^1(\R _+ ).
\end{align*}
Finally we prove \eqref{eq:asstab4}.  Like in \eqref{eq:equipart1} and using \eqref{eq:asstab3}  and \eqref{eq:expmom} we have
 \begin{align}\label{eq:equipart1}
0=\lim_{t\to +\infty}\(\mathbf{P}(u_0)-\mathbf{P}(  e^{\frac{\im}{2} v(t)x }\phi_{\omega (t)})-\mathbf{P}(\eta (t))\)= \lim_{t\to +\infty}\(\mathbf{P}(u_0)- \frac{ v(t)}{2}\mathbf{q}( {\omega _+})-\mathbf{P}(\eta (t))\)   .
\end{align}
 So the key is to show $\frac{d}{dt}\mathbf{P}(\eta) \in L^1(\R  _+)$. We have
 \begin{align*}
\frac{d}{dt}\mathbf{P}(\eta)&=-\<\im \partial _x\eta,\dot{\eta}\>
= \<\im \partial _x \eta,e^{-\im \vartheta +D\partial _x}     D _{\Theta} \phi [\Theta ] \(  \dot \Theta  -\widetilde{\Theta}   \) \>    \\&\quad -\<\im \partial _x \eta,  \im  \( f( \phi [\omega ,v , z ] + \eta  )  -  f( \phi  [\omega ,v , z ]  )  \)\> -\<\im \partial _x \eta, \im R[\omega ,v , z ]\>\\&
=J_1+J_2+J_3 .
\end{align*}
 Like for $I$ and $III$ above we have $J_1,J_3\in L^1(\R  _+)$ so we turn to  $J_2$ which we write  \small
 \begin{align*}
    J_2  &= -    \int _{0}^{1} ds\<  f( \phi [\omega ,v , z ] + \eta  )  -  f( \phi  [\omega ,v , z ])   -  Df(   s\eta  ) \eta , \partial _x\eta \>   \\&  = -\sum _{j=1,2 }   \int _{0}^{1} ds \int _{\Omega_{j,t,s }}\<  f( \phi [\omega ,v , z ] + \eta  )  -  f( \phi  [\omega ,v , z ])   -  Df(   s\eta  ) \eta , \partial _x\eta \>  _\C dx =: J_{21}+ J_{22} .
 \end{align*} \normalsize
We have
\begin{align*}
   J_{21} (t) =-\int _{0}^{1} ds  \int_{\Omega_{1,t,s }}   \< \left [ D f   ( \phi[\omega (t)  ,v(t)   ,z (t) ]+  s\eta  (t) )        -Df(   s  \eta (t) ) \right ] \eta (t) , \partial _x\eta (t)\> _\C dx.
 \end{align*}
 Then
\begin{align*}
  | J_{21} (t) |\lesssim  \int _{0}^{1} ds \int_{\Omega_{1,t,s }}   |\phi[\omega (t)  ,v(t)   ,z (t) ]   |  ^{p-1}        |\eta (t) |   | \partial _x \eta (t) |    dx  \lesssim   \|  e^{- a\<x\>}    \eta (t)\| _{H^1 (\R )}  ^2 \in L^1(\R _+ ).
 \end{align*}
 We have
 \begin{align*}
J_{22}(t)&=-  \int _{[0,1]^2} d\tau \ ds  \int_{\Omega_{2,t,s }}  \<  D ^2f( \tau \phi[\omega (t)  ,v(t)   ,z (t) ]  +s \eta (t) )   \(  \eta  (t)   , \phi[\omega (t)  ,v(t)   ,z (t) ] \)  , \partial _x\eta (t) \>  _{\C}
\end{align*} so that  we get the following, which completes the proof  of \eqref{eq:asstab4}, \begin{align*}
|J_{22}(t)|&\lesssim   \int _{[0,1]^2} d\tau \ ds \ s^{p-2} \int_{\Omega_{2,t,s }}  |  \eta (t) | ^{p-2}| \phi[\omega (t)  ,v(t)   ,z (t) ] | ^{p-1 +2-p}     \ |  \eta (t) | \   | \partial _x \eta (t)|  \  \,dx \\& \lesssim \int_{\R} |\phi[\omega (t),v(t),z(t)]|^{p-1} |\eta (t)|  | \partial _x \eta (t)| \,dx    \lesssim \|  e^{- a\<x\>}    \eta (t)\| _{H^1 (\R )}  ^2  \in L^1(\R  _+).
\end{align*}

\qed

\section{Proof of Proposition  \ref{prop:modpar} }\label{sec:modpar}

We will need the following.
\begin{lemma}\label{lem:discrestnonlin1} For $\mathbf{v} \in \{  \im  \phi [\omega ,v , z ], \partial _\omega \phi [\omega ,v , z ] ,\partial _x \phi [\omega ,v , z ]  , \im x \phi [\omega ,v , z ] \}$  we have
\begin{align}
  \label{eq:discrestnonlin1} |\<     f( \phi [\omega ,v , z ] + \eta  )  -  f( \phi  [\omega ,v , z ]  ) -D f( \phi   [\omega ,v , z ] ) \eta ,     \mathbf{v} \>|\lesssim \|  \eta  \| _{\widetilde{\Sigma}} ^2 +|z| ^4.
\end{align}
Furthermore we have
\begin{align} \label{eq:discrestnonlin1zj}
 | \<      f( \phi [\omega ,v , z ] + \eta  )  -  f( \phi  [\omega ,v , z ]  ) -D f( \phi  [\omega ,v , z ]  ) \eta ,       \partial _{z_j} \phi [\omega , v, z ]    \> |\lesssim \|  \eta  \| _{\widetilde{\Sigma}} ^ p.
\end{align}

\end{lemma}
\proof  We start with the following.
\begin{claim}\label{claim:discrestnonlin2} For $\mathbf{v} \in \{  \im  \phi [\omega ,v , z ],  \im x \phi [\omega ,v , z ] \}$  we have
\begin{align}
  \label{eq:discrestnonlin11} |\<     f( \phi [\omega ,v , z ] + \eta  )  -  f( \phi  [\omega ,v , z ]  ) -D f( \phi   [\omega ,v , z ] ) \eta ,     \mathbf{v} \>|\lesssim \|  \eta  \| _{\widetilde{\Sigma}} ^2.
\end{align}
\end{claim}
\proof  We start here with $ \mathbf{v} =  \im  \phi [\omega ,v , z ] $. Notice that the constant $\im$ does not play a particular role.
Set  for $s\in (0,1)$
\begin{align*}&
   \Omega_{1,t ,s}'=\{x\in \R\ |\  2| s \eta(t,x)| \le |\phi[\omega (t),v(t) ,z(t)]|\}  \text{ and }\\&  \Omega_{2,t,s }'=\R\setminus \Omega_{1,t,s } '=\{x\in \R\ |\ 2| s\eta(t,x)|>  |\phi[\omega (t),v(t),z(t)]|\}
\end{align*}
and  split
\begin{align*} &
  \<     f( \phi [\omega ,v , z ] + \eta  )  -  f( \phi  [\omega ,v , z ]  ) -D f( \phi   [\omega ,v , z ] ) \eta ,    \im  \phi [\omega ,v , z ] \> =I_1(t)+I_2(t) \text{ for}
   \\& I_j(t)=  \int _0^1 ds \int _{\Omega_{j,t,s }'}  \<   \left [ D f( \phi [\omega ,v , z ] + s\eta  )    -D f( \phi   [\omega ,v , z ] ) \right ] \eta ,    \im  \phi [\omega ,v , z ] \>  _\C dx .
\end{align*}
We have
\begin{align*} &
   I_1 (t)=  \int _{[0,1]^2} d \tau \  ds  \ s\int _{\Omega_{1,t,s }'}  \<   D ^2f( \phi [\omega ,v , z ] + s \tau \eta  )  ( \eta , \eta ) ,    \im  \phi [\omega ,v , z ] \>  _\C dx
\end{align*}
with
\begin{align*}
  |I_1(t)|& \lesssim    \int _{[0,1]^2} d \tau \  ds  \int _{\Omega_{1,t,s }'} | \phi [\omega ,v , z ]| ^{p-1} |\eta | ^2 d x \lesssim \|  \eta  \| _{\widetilde{\Sigma}} ^2.
\end{align*}
We have the following, which completes the proof of  \eqref{eq:discrestnonlin1}  for $ \mathbf{v} =  \im  \phi [\omega ,v , z ] $,  \small
\begin{align*}
  |I_2(t)|& \lesssim    \int _{[0,1] }  \  ds  \int _{\Omega_{2,t,s }'} | \phi [\omega ,v , z ]|   |\eta | ^p  d x \lesssim    \int _{[0,1] }  \  ds  \int _{\Omega_{2,t,s }'} | \phi [\omega ,v , z ]| ^{p-1}   |\eta | ^{2-p}  |\eta | ^p  d x   \le
   \|  \eta  \| _{\widetilde{\Sigma}} ^2  .
\end{align*}
\normalsize
For $ \mathbf{v} =  \im x \phi [\omega ,v , z ] $ the argument is similar.

\qed

We continue the proof of Lemma \ref{lem:discrestnonlin1} considering   the case $\mathbf{v}=\partial _x \phi [\omega ,v , z ] = \frac{\im}{2}v   \phi [\omega ,v , z ] + e^{\frac{\im}{2} vx}   \partial _x       \phi [\omega   , z ] $.
It is enough to focus on the following terms
\begin{align*} &
 \int _0^1 ds  \<      \left [ D f( \phi [\omega ,v , z ] + s\eta  )    -D f( \phi   [\omega ,v , z ] ) \right ] \eta , e^{\frac{\im}{2} vx}   \partial _x       \phi [\omega   , z ] \> _{L^2\( \Omega_{1,t  } '' \) } \\& +  \int _0^1 ds  \<      \left [ D f( \phi [\omega ,v , z ] + s\eta  )    -D f( \phi   [\omega ,v , z ] ) \right ] \eta ,  e^{\frac{\im}{2} vx}      \partial _x  \phi [\omega   , z ] \> _{L^2\( \Omega_{2,t  } '' \) }  :=K_1+K_2
\end{align*}
in corresponding to a   partition  \begin{align*}&
   \Omega_{1,t  }''=\{x\in \R\ |\  |   \phi _\omega |\gg |   \widetilde{\phi} [\omega   , z ]|  \text{ and }   | \partial _x \phi _\omega |\gg | \partial _x  \widetilde{\phi} [\omega   , z ]| \}  \text{ and } \Omega_{2,t  }''=\R \backslash  \Omega_{1,t  }'' .
\end{align*}
   Since
$  | \partial _x \phi _\omega | \lesssim | \phi _\omega |$  for $K_1$ we can repeat the arguments    in  Claim \ref{claim:discrestnonlin2} replacing in the inequalities $\partial _x  \phi [\omega   , z ]$   with  $  \phi [\omega   , z ]$.  So $K_1$ satisfies the same bounds in \eqref{eq:discrestnonlin11}.  Next, we split
\begin{align*}
  K_2= K _{3}+K _{4}
\end{align*}
in correspondence to  a   partition
\begin{align*}&
    \Omega_{2,t  }'' =  \Omega_{3,t  }'' \cup \Omega_{4,t  }''   \text{   with } \\&  \Omega_{3,t  }''=\{ x\in  \Omega_{2,t  }'' \ |\   \phi _\omega |\lesssim  |   \widetilde{\phi} [\omega   , z ]|  \text{ and }   | \partial _x \phi _\omega |\gg  | \partial _x  \widetilde{\phi} [\omega   , z ]| \}   \text{   and } \\&  \Omega_{4,t  }''=\{ x\in  \Omega_{2,t  }'' \ |\       | \partial _x \phi _\omega |\lesssim  | \partial _x  \widetilde{\phi} [\omega   , z ]| \} .
\end{align*}
Then for $K_3$ we can bound $   |\partial _x  \phi [\omega   , z ]|\lesssim  | \partial  _ x \phi _\omega |\lesssim  |  \phi _\omega | \lesssim  |   \widetilde{\phi} [\omega   , z ]|$  and we have
\begin{align}\label{eq:K3est}
  |K_3|\le  \int _0^1 ds \int _{\Omega_{3,t  }''} \left |   D f( \phi [\omega ,v , z ] + s\eta  )    -D f( \phi   [\omega ,v , z ] )  \right | \ |\eta | \ |     \widetilde{\phi} [\omega   , z ]| dx.
\end{align}

\begin{claim}
  \label{claim:estK3}  $|K_3|\lesssim \|  \eta  \| _{\widetilde{\Sigma}} ^2 +|z| ^4$.
\end{claim}
\proof We split
\begin{align*}
  \Omega_{1,3,t ,s }''=\{   x\in \Omega_{3,t  }'': \quad |s \eta | \gtrsim |   {\phi} [\omega   , z ]| \} \text{ and }\Omega_{2,3,t ,s }''=\{   x\in \Omega_{3,t  }'': \quad |s \eta | \ll |   {\phi} [\omega   , z ]| \} .
\end{align*}
Then we split
$K_{3}=K_{31}+K_{32}$
\begin{align*}
 |K_{31}|\lesssim  \int _{[0,1] }   ds \int _{ \Omega_{1,3,t ,s }''}  |\eta | ^p   |     \widetilde{\phi} [\omega   , z ]|    ds  \lesssim \int _{   \R} e^{-|x|\sqrt{\omega} \sqrt{1-\lambda (p)}}  \( |z|^4 + |\eta | ^{p \frac{4}{3}} \)  \lesssim  \|  \eta  \| _{\widetilde{\Sigma}} ^2 +|z| ^4
\end{align*}
since we have $p >3/2$. Next we have
\begin{align*}
    |K_{32}|&\lesssim  \int _{[0,1] ^2 }   ds \ d\tau  \int _{ \Omega_{1,3,t ,s }''} \left |   D ^2 f( \phi [\omega ,v , z ] + \tau s\eta  )     \right | \ |s\eta |    \ |\eta | \ |     \widetilde{\phi} [\omega   , z ]| dx \\& \lesssim \int _{[0,1] ^2 }   ds \ d\tau  \int _{ \Omega_{1,3,t ,s }''} \left |    \phi [\omega   , z ]   \right | ^{p-2} \ |s\eta | ^{p-1+2-p}   \ |\eta | \ |     \widetilde{\phi} [\omega   , z ]| dx \\& \lesssim   \int _{  \R}   |\eta | ^p   |     \widetilde{\phi} [\omega   , z ]|  \lesssim  \|  \eta  \| _{\widetilde{\Sigma}} ^2 +|z| ^4.
\end{align*}
\qed

 Now we consider $K_4$. We have
\begin{align*}
  |K_4|\le  \int _0^1 ds \int _{\Omega_{4,t  }''} \left |   D f( \phi [\omega ,v , z ] + s\eta  )    -D f( \phi   [\omega ,v , z ] )  \right | \ |\eta | \ |  \partial _x  \widetilde{\phi} [\omega   , z ]| dx .
\end{align*}
Splitting
\begin{align*}
  \Omega_{1,4,t ,s }''=\{   x\in \Omega_{4,t  }'': \quad |s \eta | \gtrsim |   {\phi} [\omega   , z ]| \} \text{ and }\Omega_{2,4,t ,s }''=\{   x\in \Omega_{4,t  }'': \quad |s \eta | \ll |   {\phi} [\omega   , z ]| \}
\end{align*}
then the proof of  $|K_4|\lesssim  \|  \eta  \| _{\widetilde{\Sigma}} ^2 +|z| ^4$ is similar to the proof of Claim \ref{claim:estK3}.

Next, continuing the proof of  proof of Lemma \ref{lem:discrestnonlin1} we   consider case
\begin{align*}
 \mathbf{v}&=   e^{\frac{\im}{2} v x}  \partial _\omega  \(   \omega^{\frac{1}{p-1}} \phi [  z ] (\sqrt{\omega}x)  \) =  e^{\frac{\im}{2} v x}  \omega ^{\frac{2-p}{p-1}}  (\Lambda_p    \phi [  z ])(\sqrt{\omega}x) .
\end{align*}
 By the definition of $\Lambda_p $ in \eqref{eq:lambdap} we are reduced to the previous cases.

 \noindent Finally, we look at the proof of \eqref{eq:discrestnonlin1zj}. For $\mathbf{v}=\partial _{z_j} \phi [\omega , v, 0 ]$   We need to bound
 \begin{align*}
  \sum _{j=1 ,2} \int _0 ^1 ds  \<     \( D f( \phi [\omega ,v , z ] + s\eta  )    -D f( \phi  [\omega ,v , z ]  )\) \eta ,       \mathbf{v}    \> _{L^2\(  \Omega_{j,t ,s}'  \)} =:L_1+L_2.
 \end{align*}
 Then
 \begin{align*}
   |L_1(t)|&\lesssim  \int _{[0,1]^2} d \tau \  ds  \ s\int _{\Omega_{1,t,s }'} | \<   D ^2f( \phi [\omega ,v , z ] + s \tau \eta  )  ( \eta , \eta ) ,    \mathbf{v} \>  _\C | dx\\& \lesssim  \int _{[0,1]^2} d \tau \  ds  \ s\int _{\Omega_{1,t,s }'} | \phi [\omega   , z ]| ^{p-2} |s\eta| ^{2-p+p-1}  \  | \eta|  e^{-|x| \sqrt{\omega} \sqrt{1-\lambda (p)}}      dx \\& \lesssim \int _\R  | \eta| ^{p}  e^{-|x| \sqrt{\omega} \sqrt{1-\lambda (p)}}      dx \lesssim \| \eta \| _{\widetilde{\Sigma}} ^p.
 \end{align*}
Finally we have the following, which completes the proof of \eqref{eq:discrestnonlin1zj},
 \begin{align*}
   |L_2(t)|&\lesssim  \int _{[0,1] }    ds   \int _{\Omega_{2,t,s }'} |s\eta | ^{p-1}  \ |\eta | \    |    \mathbf{v}  | dx  \lesssim \| \eta \| _{\widetilde{\Sigma}} ^p.
 \end{align*}

\qed

\begin{lemma}\label{lem:lemdscrt} We have the estimates
 \begin{align} \label{eq:discrest1}
   &|\dot \vartheta -\widetilde{\vartheta} | + |\dot \omega -\widetilde{\omega} | + |\dot D -v | +  |\dot v   |  \lesssim   |z|^4+    \|   \eta \| _{\widetilde{\Sigma}}     ^2  +\| \eta \| _{\Sigma_{A}}^2 ,   \\&  |\dot z -\widetilde{z} | \lesssim   |z| \|   \eta \| _{\widetilde{\Sigma}}     +   \|   \eta \| _{\widetilde{\Sigma}}  ^p +\| \eta \| _{\Sigma_{A}}^2+ |z| ^4  .   \label{eq:discrest2}
\end{align}

\end{lemma}
\proof
Applying $\< \cdot ,\im e^{-\im \vartheta +D\partial _x}  D_{\Theta} \phi [\Theta ]   \boldsymbol{\Theta } \>  $  with $\boldsymbol{\Theta }\in \R ^4\times \C$   to \eqref{eq:nls3}
we have
\begin{align}\label{eq:discrest01l1}   &  \<   \dot \eta  , \im e^{-\im \vartheta +D\partial _x}  D_{\Theta} \phi [\Theta ]   \boldsymbol{\Theta } \>  +    \<    \( \im  \dot \vartheta  -  \dot D    \partial _x \) \eta  , \im e^{-\im \vartheta +D\partial _x}  D_{\Theta} \phi [\Theta ]   \boldsymbol{\Theta } \>     \\& +         \< D _{\Theta} \phi [\Theta ] (\dot \Theta -\widetilde{\Theta} ) , \im  D _{\Theta} \phi [\Theta ]   \boldsymbol{\Theta } \>
      \label{eq:discrest01l2} \\& =   \< \eta ,   \(\partial _x^2    +D f( \phi [\omega ,v , z ]  )  \)  e^{-\im \vartheta +D\partial _x}   D _{\Theta} \phi [\Theta ]   \boldsymbol{\Theta } \>  \nonumber \\&   -\<     f( \phi [\omega ,v , z ] + \eta  )  -  f( \phi  [\omega ,v , z ]  ) -D f( \phi  [\omega ,v , z ]  ) \eta ,    e^{-\im \vartheta +D\partial _x}   D _{\Theta} \phi [\Theta ]   \boldsymbol{\Theta } \>    \label{eq:discrest01l3}
   \end{align}
where the contribution of $R[\omega ,v   , z ]$ cancels by \eqref{R:orthboost}.
Now we have
\begin{align}  \label{eq:discrest01l4}
    \<   \dot \eta  , \im e^{-\im \vartheta +D\partial _x}  D_{\Theta} \phi [\Theta ]   \boldsymbol{\Theta } \> &= -\<     \eta  , \im e^{-\im \vartheta +D\partial _x}  D_{\Theta}^2 \phi [\Theta ]   \(\boldsymbol{\Theta } , \dot \Theta \) \>\\&  +  \<    \eta  , \im      \( \im  \dot \vartheta  -  \dot D    \partial _x \) e^{-\im \vartheta +D\partial _x}  D_{\Theta} \phi [\Theta ]   \boldsymbol{\Theta } \> .\label{eq:discrest01l5}
\end{align}
Setting also  $R[\Theta]  =  e^{ \im \vartheta -D\partial _x} R[\omega,v, {z}]$, notice that equation  \eqref{eq:refboost} can be written as
\begin{align*}
  \partial ^2_x\phi  [\Theta]+ f(\phi [\Theta]) +\im  D _{\Theta}   \phi [\Theta ]  \widetilde{\Theta}  =R[\Theta] .
\end{align*}
Differentiating  in $\Theta$,   we obtain
\begin{align*}
  \( \partial ^2_x   + Df(\phi [\Theta]) \) D _{\Theta} \phi [\Theta ]   \boldsymbol{\Theta } +\im  D _{\Theta}^2   \phi [\Theta ]  (\widetilde{\Theta}  , \boldsymbol{\Theta })+  \im  D _{\Theta}   \phi [\Theta ]  D _{\Theta}  \widetilde{\Theta} \boldsymbol{\Theta }   =D _{\Theta} R[\Theta] \boldsymbol{\Theta }.
\end{align*}
Since   $f( \phi [\Theta]  )  = e^{ \im \vartheta - D\partial _x } f(   \phi [\omega ,v , z ] )$ yields $e^{ -\im \vartheta + D\partial _x } D f( \phi [\Theta] ) X=  D f(  \phi [\omega ,v , z ] ) e^{ -\im \vartheta + D\partial _x }  X$,  we get  \small
\begin{align*}&
      \<     \eta,  \(\partial _x^2    +D f( \phi [\omega ,v , z ] )  \)  e^{ -\im \vartheta + D\partial _x }  D _{\Theta} \phi [\Theta ]   \boldsymbol{\Theta } \>          =   \<  e^{  \im \vartheta - D\partial _x }\eta ,  \(  \partial ^2_x + Df(\phi [\Theta])  \) D _{\Theta} \phi [\Theta ]   \boldsymbol{\Theta } \>   \\& = -\<  e^{  \im \vartheta - D\partial _x } \eta , \im  D _{\Theta}^2   \phi [\Theta ]  (\widetilde{\Theta}  , \boldsymbol{\Theta })+ \cancel{ \im  D _{\Theta}   \phi [\Theta ]  D _{\Theta}  \widetilde{\Theta} \boldsymbol{\Theta }} -D _{\Theta} R[\Theta] \boldsymbol{\Theta }\>
\end{align*}
\normalsize  where the cancellation follows by the modulation orthogonality \eqref{61}. Entering this and \eqref{eq:discrest01l4}  in \eqref{eq:discrest01l1}--\eqref{eq:discrest01l3}
and with the cancellation of \eqref{eq:discrest01l5}  with the second term in line \eqref{eq:discrest01l1}
we obtain
\small
\begin{align}\label {eq:discrest2}   &\< D _{\Theta} \phi [\Theta ] (\dot \Theta -\widetilde{\Theta} ) , \im  D _{\Theta} \phi [\Theta ]   \boldsymbol{\Theta } \>
  -      \<\eta , \im   e^{-\im \vartheta +D\partial _x}       D _{\Theta}^2\phi [\Theta ]  ( \boldsymbol{\Theta } ,   \dot \Theta -\widetilde{\Theta}  ) \>    \\& \nonumber \cancel{ -      \<\eta , \im   e^{-\im \vartheta +D\partial _x}       D _{\Theta}^2\phi [\Theta ]  ( \boldsymbol{\Theta } ,    \widetilde{\Theta}  ) \> }     = \cancel{ -\<  e^{  \im \vartheta - D\partial _x } \eta , \im  D _{\Theta}^2   \phi [\Theta ]  (\widetilde{\Theta}  , \boldsymbol{\Theta }) \>  } +  \<  e^{  \im \vartheta - D\partial _x } \eta ,   D _{\Theta} R[\Theta] \boldsymbol{\Theta } \> \nonumber \\&   -\<     f( \phi [\omega ,v , z ] + \eta  )  -  f( \phi  [\omega ,v , z ]  ) -D f( \phi  [\omega ,v , z ]  ) \eta ,    e^{-\im \vartheta +D\partial _x}   D _{\Theta} \phi [\Theta ]   \boldsymbol{\Theta } \>   \nonumber
   \end{align}
where the two cancelled terms cancel each other. From this  we get
\begin{align*}   &
    -  (\dot \omega -\widetilde{\omega}   )  \< \partial _\omega \phi [\omega ,v , z ] , \phi  [\omega ,v , z ] \> +  (\dot \vartheta -\widetilde{\vartheta}   )  \cancel{\<   \phi [\omega ,v , z ] , \im \phi  [\omega ,v , z ] \>}
   \\&   +(\dot D -\widetilde{D}   )\cancel{ \< \partial _x \phi [\omega ,v , z ] , \phi  [\omega ,v , z ] \>}
    -  (\dot v -\widetilde{v}   ) \cancel{ \<  \im  \frac{x}{2}\phi [\omega ,v , z ] ,   \phi  [\omega ,v , z ] \>}
    \\&- \< D _{z} \phi [\omega ,v , z ]  (\dot z -\widetilde{z})  , \phi [\omega ,v , z ] \> + O\(   \|   \eta \| _{\widetilde{\Sigma}}       |\dot \Theta -\widetilde{\Theta} |\) \\& = \<   \eta ,
  \im  R [\omega ,v , z ]   \>        -\<     f( \phi [\omega ,v , z ] + \eta  )  -  f( \phi  [\omega ,v , z ]  ) -D f( \phi   [\omega ,v , z ] ) \eta ,    \im  \phi [\omega ,v , z ] \>
\end{align*}
where the cancelled terms are null for elementary reasons  related also for various cancellations in the sequel of this lemma also to the fact that by \eqref{eq:sol} and Lemma \ref{lem:simandmult} the function $ \phi [  {z}]$ is even in $x$. From  \eqref{estR} and \eqref{eq:discrestnonlin1} we obtain
\begin{align} \label{eq:omega}  &
      (\dot \omega -\widetilde{\omega}   )     q' (\omega )  + O\(   \( | z|   +   \|   \eta \| _{\widetilde{\Sigma}} \)       |\dot \Theta -\widetilde{\Theta} |\)  =   O\(   \|   \eta \| _{\widetilde{\Sigma}}    |z|^2\)  +  O\(  \|   \eta \| _{\widetilde{\Sigma}}    ^2\)  .
\end{align}
From  \eqref {eq:discrest2}   we obtain
\begin{align*}   &
      (\dot \vartheta  -\widetilde{\vartheta}   )  \<   \phi [\omega ,v , z ] , \partial _\omega  \phi  [\omega ,v , z ] \> +  (\dot \omega  -\widetilde{\omega}   )  \cancel{\<  \partial _\omega \phi [\omega ,v , z ] , \im  \partial _\omega\phi  [\omega ,v , z ] \>}
   \\&  -  (\dot D -\widetilde{D}   )  \< \partial _x \phi [\omega ,v , z ] , \im  \partial _\omega \phi  [\omega ,v , z ] \>
    +  (\dot v -\widetilde{v}   ) \bcancel{ \<  \im  \frac{x}{2}\phi [\omega ,v , z ] , \im  \partial _\omega\phi  [\omega ,v , z ] \>}
    \\&+ \< D _{z} \phi [\omega ,v , z ]  (\dot z -\widetilde{z})  ,\im  \partial _\omega  \phi [\omega ,v , z ] \> + O\(   \|   \eta \| _{\widetilde{\Sigma}}       |\dot \Theta -\widetilde{\Theta} |\) \\& = \<   \eta ,
     \partial _\omega  R [\omega ,v , z ]   \>        -\<     f( \phi [\omega ,v , z ] + \eta  )  -  f( \phi  [\omega ,v , z ]  ) -D f( \phi   [\omega ,v , z ] ) \eta ,       \partial _\omega   \phi [\omega ,v , z ] \>
\end{align*}
where the  canceled terms are zero by elementary reasons. We have
\begin{align*}&
  \< \partial _x \phi [\omega ,v , z ] , \im  \partial _\omega \phi  [\omega ,v , z ] \> =   \< \partial _x  \( e^{\frac{\im}{2} vx }   \phi [\omega  , z ]  \), \im  e^{\frac{\im}{2} vx } \partial _\omega \phi  [\omega   , z ] \>  \\& = \cancel{ \< \partial _x    \phi [\omega  , z ]   , \im    \partial _\omega \phi  [\omega   , z ] \>} + \frac{v}{2}  \<      \phi [\omega  , z ]   ,    \partial _\omega \phi  [\omega   , z ] \>
  =  \frac{v}{2}  \( \mathbf{q}' (\omega ) +O(|z|) \)  ,
\end{align*}
where the canceled term is null.

We have
\begin{align*} &
  \<  \eta ,
     \partial _\omega  R [\omega ,v , z ]   \>  =  \omega ^{\frac{1}{p-1}}   \<   e^{\frac{\im}{2} vx } \eta , \left [ \(  \frac{p}{p-1} + \frac{x}{2}\partial _x\)   R [  z ] \right ]  \(  \sqrt{\omega} x  \) \> \\& = \omega ^{\frac{1}{p-1} + \frac{1}{2}} \<    \(  \frac{p}{p-1} + \frac{x}{2}\partial _x\)  ^* \(  \(  e^{\frac{\im}{2} vx } \eta  \) \( \omega ^{-\frac{1}{2}} x\)\) ,    R [  z ]  \>
\end{align*}
so that it is elementary to conclude that
\begin{align*}&
  |\<  \eta ,
     \partial _\omega  R [\omega ,v , z ]   \> | \lesssim \( \| \eta \| _{\Sigma_{A}} + \| \eta \| _{\widetilde{\Sigma} } \) |z|^2.
\end{align*}
We obtain
\begin{align}  \nonumber  &
     \( (\dot \vartheta -\widetilde{\vartheta}   )   - \frac{v}{2}  (\dot D -\widetilde{D}   )  \)   \mathbf{q}' (\omega )    +  O\(   \( | z|   +   \|   \eta \| _{\widetilde{\Sigma}} \)       |\dot \Theta -\widetilde{\Theta} |\) \\& =   O\(   \|   \eta \| _{\widetilde{\Sigma}} ^2 +   | z|^4\)  . \label{eq:theta}
\end{align}
From \eqref{eq:discrest2}  we have
\begin{align*}  &  (\dot \vartheta  -\widetilde{\vartheta}   ) \cancel { \<   \phi [\omega ,v , z ] ,\partial _x \phi [\omega ,v , z ]   \> }+  (\dot \omega  -\widetilde{\omega}   )  \<  \partial _\omega \phi [\omega ,v , z ] , \im  \partial _x \phi [\omega ,v , z ]  \>
   \\&  -  (\dot D -\widetilde{D}   ) \cancel{ \< \partial _x \phi [\omega ,v , z ] , \im \partial _x \phi [\omega ,v , z ]  \> }
    -  (\dot v -\widetilde{v}   )   \( \frac{\mathbf{q}(\omega )}{2} + O(z) \)
    \\&+ \< D _{z} \phi [\omega ,v , z ]  (\dot z -\widetilde{z})  ,\im \partial _x \phi [\omega ,v , z ]   \>          +O\(   \|   \eta \| _{\widetilde{\Sigma}}        |\dot \Theta -\widetilde{\Theta} |\)\\& \nonumber     =
\<     \eta ,
  \partial _{x} R [\omega , v, z ]  \>       -\<      f( \phi [\omega ,v , z ] + \eta  )  -  f( \phi  [\omega ,v , z ]  ) -D f( \phi  [\omega ,v , z ]  ) \eta ,       \partial _{x} \phi [\omega , v, z ]    \>,   \nonumber
   \end{align*}
from which we derive
\begin{align}  \nonumber  &  -\frac{v}{2} (\dot \omega  -\widetilde{\omega}   )   \mathbf{q}' (\omega )
  -  (\dot v -\widetilde{v}   )   \frac{\mathbf{q}(\omega )}{2}
\\&   +  O\(    \( | z|   +  \|   \eta \| _{\widetilde{\Sigma}}   \)     |\dot \Theta -\widetilde{\Theta} |\)   =   O\(  \|   \eta \| _{ {\Sigma}_A} ^2 +  \|   \eta \| _{\widetilde{\Sigma}} ^2 +   | z|^4\)  . \label{eq:eqv}
\end{align}
From \eqref{eq:discrest2}  we have
\begin{align*}  & - (\dot \vartheta  -\widetilde{\vartheta}   )  \cancel{ \< \im   \phi [\omega ,v , z ] , x \phi [\omega ,v , z ]   \> }-  (\dot \omega  -\widetilde{\omega}   ) \cancel{ \<  \partial _\omega \phi [\omega ,v , z ] , x \phi [\omega ,v , z ]  \>}
   \\&  +  (\dot D -\widetilde{D}   )   \<  \partial _x \phi [\omega ,v , z ] , x   \phi [\omega ,v , z ]  \>
    -  (\dot v -\widetilde{v}   )    \cancel{ \< \frac{\im}{2} x\phi [\omega ,v , z ] , x   \phi [\omega ,v , z ]  \> }
    \\&- \< D _{z} \phi [\omega ,v , z ]  (\dot z -\widetilde{z})  , x \phi [\omega ,v , z ]   \>          +O\(   \|   \eta \| _{\widetilde{\Sigma}}        |\dot \Theta -\widetilde{\Theta} |\)\\& \nonumber     =
\<     \eta ,
  \partial _{v} R [\omega , v, z ]  \>       -\<      f( \phi [\omega ,v , z ] + \eta  )  -  f( \phi  [\omega ,v , z ]  ) -D f( \phi  [\omega ,v , z ]  ) \eta ,         {x} \phi [\omega , v, z ]    \>     \nonumber
   \end{align*}
where the cancelled terms are null for elementary reasons. This yields
\begin{align}
  \label{eq:eqD} (\dot D -\widetilde{D}   ) \frac{q(\omega )}{2}  =  O\(    \( | z|   +  \|   \eta \| _{\widetilde{\Sigma}}   \)     |\dot \Theta -\widetilde{\Theta} |\) +  O\(  \|   \eta \| _{ {\Sigma}_A} ^2 +  \|   \eta \| _{\widetilde{\Sigma}} ^2 +   | z|^4\)
\end{align}
where we used
\begin{align*} &
   \<     \eta ,
  \partial _{v} R [\omega , v, z ]  \> =  \<     \eta , \frac{\im}{2}x e^{\frac{\im}{2}xv}
    R [\omega ,   z ]  \> = O\(   \|   \eta \| _{\widetilde{\Sigma}}     | z|^2 \)
\end{align*}
and Lemma \ref{lem:discrestnonlin1},
From \eqref{eq:discrest2}  we have
\begin{align*}  &  \< D _{z} \phi [\omega ,v , z ]  (\dot z -\widetilde{z})  ,\im  \partial _{z_j} \phi [\omega ,v , z ] \>          +O\( \( |z| +  \|   \eta \| _{\widetilde{\Sigma}}\)         |\dot \Theta -\widetilde{\Theta} |\)\\& \nonumber     =
\<     \eta ,
  \partial _{z_j} R [\omega , v, z ]  \>   \\&    -\<      f( \phi [\omega ,v , z ] + \eta  )  -  f( \phi  [\omega ,v , z ]  ) -D f( \phi  [\omega ,v , z ]  ) \eta ,       \partial _{z_j} \phi [\omega , v, z ]    \> .   \nonumber
   \end{align*}
\begin{claim}
  \label{claim:estdzjR} We have
  \begin{align}
    \label{eq:estdzjR} | \<     \eta ,\partial _{z_j} R [\omega , v, z ]  \>|\lesssim \| \eta \| _{\widetilde{\Sigma}}  \  |z| .
  \end{align}
\end{claim}
\proof   From \eqref{eq:refboost} we can write  \small
\begin{align*}
  R[ \omega , v,  {z}]=\widehat{R}[ \omega , v,  {z}] - \widetilde{\vartheta}_{\mathcal{R}}[ \omega ,v,  {z}]\phi [ \omega  ,v , {z}]+ \im \widetilde{\omega}_{\mathcal{R}}[ \omega ,v,  {z}] \partial _\omega \phi [  \omega ,v, {z}]   +\im  D_{z}\phi  [ \omega , v, {z}]\widetilde{z }_{\mathcal{R}}[ \omega ,v,  {z}]
\end{align*} \normalsize
where   $\widehat{R}[ \omega , v,  {z}] = e^{\frac{\im}{2} vx}\omega ^{\frac{p}{p-1}} \widehat{R}[     {z}]  \( \sqrt{\psi} x\)$.
Notice that by elementary computations we have
\begin{align*}
 \widehat{R}[ \omega , v,  {z}]=f(  \phi [\omega , v] +   \widetilde{ \phi} [    \omega , v, {z}]) -  f(\phi [    \omega , v ])          - Df(\phi  [    \omega , v ] )\widetilde{\phi}    [    \omega , v, {z}]
\end{align*}
where $\phi [\omega , v] := e^{\frac{\im}{2} vx} \phi _\omega$.
 By Proposition \ref{prop:refpropf} we have  \small
\begin{align*} &
 | \<     \eta ,\partial _{z_j} \(  - \widetilde{\vartheta}_{\mathcal{R}}[ \omega ,v,  {z}]\phi [ \omega  ,v , {z}]+ \im \widetilde{\omega}_{\mathcal{R}}[ \omega ,v,  {z}] \partial _\omega \phi [  \omega ,v, {z}]   +\im  D_{z}\phi  [ \omega , v, {z}]\widetilde{z }_{\mathcal{R}}[ \omega ,v,  {z}]\)  \> |   \lesssim \| \eta \| _{\widetilde{\Sigma}}  \  |z| .
\end{align*}
\normalsize
Finally, we have
\begin{align*} &
 | \<     \eta ,\partial _{z_j} \widehat{R} [\omega , v, z ]  \> |=   |\<     \eta ,\partial _{z_j}  \( f(  \phi [\omega , v] +   \widetilde{ \phi} [    \omega , v, {z}]) -  f(\phi [    \omega , v ])          - Df(\phi  [    \omega , v ] )\widetilde{\phi}    [    \omega , v, {z}]   \)  \> | \\& =|\<     \eta , \( D f(  \phi [\omega , v] +   \widetilde{ \phi} [    \omega , v, {z}])   - Df(\phi  [    \omega , v ] ) \partial _{z_j}  \) \widetilde{\phi}  [     \omega , v, 0 ]    \> |\\& \lesssim \int _{\R} |\eta |  |z|   e^{|x| \( 2-p-2\sqrt{1-\lambda (p)} \)  } \lesssim \| \eta \| _{\widetilde{\Sigma}}  \  |z| .
\end{align*}
 \qed

 From Claim \ref{claim:estdzjR} and from \eqref{eq:discrestnonlin1zj} we derive
\begin{align*}
   &  \< D _{z} \phi [\omega ,v , z ]  (\dot z -\widetilde{z})  ,\im  \partial _{z_j} \phi [\omega ,v , z ] \>          +O\( \( |z| +  \|   \eta \| _{\widetilde{\Sigma}}\)         |\dot \Theta -\widetilde{\Theta} |\)\\& \nonumber     =
  O\(  \| \eta \| _{\widetilde{\Sigma}}  \  |z|  \) +  O\(  \| \eta \| _{\widetilde{\Sigma}}  ^p \) .   \nonumber
\end{align*}
From \eqref{eq:matrixa1} the above implies
\begin{align}  & |\dot z -\widetilde{z} |\lesssim           O\( \( |z| +  \|   \eta \| _{\widetilde{\Sigma}}\)         |\dot \Theta -\widetilde{\Theta} |\)   +
  O\(  \| \eta \| _{\widetilde{\Sigma}}  \  |z|  \) +  O\(  \| \eta \| _{\widetilde{\Sigma}}  ^p \)     \label{eq:ineqz0}
   \end{align}
and so also
\begin{align} \label{eq:ineqz} &  |z| \ |\dot z -\widetilde{z} |\lesssim           O\( |z|  \( |z| +  \|   \eta \| _{\widetilde{\Sigma}}\)         |\dot \Theta -\widetilde{\Theta} |\)   +
  O\(  \| \eta \| _{\widetilde{\Sigma}}  ^2+   |z| ^4  \)  .
   \end{align}
Then from  \eqref{eq:theta}, \eqref{eq:omega}, \eqref{eq:eqv}, \eqref{eq:eqD}  and    \eqref{eq:ineqz} we derive
\begin{align*}
  |\dot \vartheta -\widetilde{\vartheta} | + |\dot \omega -\widetilde{\omega} | + |\dot D -\widetilde{D} | +  |\dot v -\widetilde{v} | + |z| \ |\dot z -\widetilde{z} | \lesssim   |z|^4+    \|   \eta \| _{\widetilde{\Sigma}}     ^2  +\| \eta \| _{\Sigma_{A}}^2
\end{align*}
which yields \eqref{eq:discrest1} and entered in \eqref{eq:ineqz0} yields \eqref{eq:discrest2}.

 \qed

\textit{Proof of Proposition  \ref{prop:modpar}.}
Lemma \ref{lem:lemdscrt} and \eqref{eq:main2} imply immediately \eqref{eq:modpar1}--\eqref{eq:modpar2}. Entering this,  \eqref{eq:estpar}  and \eqref{eq:oertstab1} in \eqref{eq:ineqz0}  we obtain \eqref{eq:modpar3}.

\section{The Fermi Golden Rule: proof of Proposition \ref{prop:FGR}}\label{sec:fgr}

  Like in Kowalczyk and Martel \cite{KM22} and Martel \cite{Martel2} and exactly like  in \cite{CM24D1} we consider the  functional
\begin{align}\label{eq:FGRfunctional}
\mathcal{J}_{\mathrm{FGR}}:=\< J   {\eta},\chi_A \(  {z}^{2} g ^{(\omega)}+ \overline{{z}}^{2} \overline{g}^{(\omega)} \)      \>  ,
\end{align}
with a nonzero  $g ^{(\omega)}\in L^\infty(\R , \C ^2)$  satisfying
\begin{align}
  \label{eq:eqsatg2} \mathcal{L}_{\omega}g ^{(\omega)}=2\im \lambda (p, \omega) g ^{(\omega)}
\end{align}
where $\lambda (p, \omega) = \omega \lambda (p) $. Here $g ^{(\omega)} (x):= g\( \sqrt{\omega}x \)$ where $g   $  solves \eqref{eq:eqsatg2}  for $\omega=1$.
That $g^{(\omega)} $ exists is known since   Krieger and Schlag \cite{KrSch}, see  Lemma 6.3, or earlier Buslaev and Perelman \cite{BP1}.

\noindent We define the FGR constant     $\gamma( p)$,   for $g   =\( g_1   ,g_2     \) ^\intercal $,   by
\begin{equation}\label{eq:fgrgamma}
   \gamma ( p )  :=     \<    \phi   ^{p-2}  \( p  \xi _1  ^2 +  \xi _2  ^2  \) , g_1  \> +2 \<    \phi   ^{p-2}    \xi _1   \xi _2  , g_2   \>.
\end{equation}
Here notice that we will prove Proposition \ref{prop:FGR} later in    section \ref{sec:prop:FGR}.
    In the next lemma we will need the following
reformulation of equation \eqref{eq:nls3}, where we identify $J=-\im $,
\small \begin{align}\label{eq:nls4}
   \dot \eta   &=   \mathcal{L}_{\omega}\eta  -  J  ( \widetilde{\vartheta}  _{\mathcal{R}} + \widetilde{\vartheta} - \dot \vartheta) \eta +   \dot D    \partial _x\eta  -   e^{J \vartheta +D\partial _x} D_\Theta \phi [\Theta]  (\dot \Theta  -\widetilde{\Theta})   +J  \( D f( \phi  [\omega ,v , z ]  ) - D f( \phi   _{ \omega }   )\) \eta
    \\&-J \( f( \phi [\omega ,v , z ] + \eta  )  -  f( \phi  [\omega ,v , z ]  ) -D f( \phi  [\omega ,v , z ]  ) \eta\) -J R[\omega ,v , z ].  \nonumber
\end{align}\normalsize
We have  the following.
\begin{lemma}\label{lem:FGR1}
There is a  $C^1 $  in time    function $\mathcal{I}_{\mathrm{FGR}}$, which satisfies $|\mathcal{I}_{\mathrm{FGR}}|\lesssim  \sqrt{\delta}$
such that
\begin{align}
\left|\dot{\mathcal{J}}_{\mathrm{FGR}}  + \dot{\mathcal{I}}_{\mathrm{FGR}}  -\omega ^{\frac{p+1}{2( p-1)}} \gamma ( p  ) |z|^4     \right |    \lesssim  A ^{-1}   \(  |z|^4   + \| \eta \| _{\Sigma _A}^{2}+ \| \eta \| ^{2}_{\widetilde{\Sigma}}  \)    . \label{eq:lem:FGR11}
\end{align}
\end{lemma}
\proof
Differentiating $\mathcal{J}_{\mathrm{FGR}}$, we have \small
\begin{align*}
\dot{\mathcal{J}}_{\mathrm{FGR}}=&
\<  J  \dot{  \eta },\chi_A \(  {z}^{2} g ^{(\omega)}+ \overline{ z }^{2} \overline{g} ^{(\omega)}\)      \>
+\< J    \eta , \chi_A \( 2z    \widetilde{z}  g^{(\omega)} +  2\overline{z}   \overline{\widetilde{z}}  \overline{g} ^{(\omega)}\)            \>   \\&+ \<  J    \eta , \chi_A \( 2z  (\dot{z}  - \widetilde{z}) g ^{(\omega)}+  2\overline{z}  (\dot{\overline{z}}  - \overline{\widetilde{z}}) \overline{g}^{(\omega)} \)            \> \\&  + \<  J    \eta , \chi_A \({z}^{2} \partial _\omega g ^{(\omega)}+ \overline{ z }^{2}\partial _\omega \overline{g} ^{(\omega)}\)      \> (\dot \omega -\widetilde{\omega}) +    \<  J    \eta , \chi_A \({z}^{2} \partial _\omega g ^{(\omega)}+ \overline{ z }^{2}\partial _\omega \overline{g} ^{(\omega)}\)      \>  \widetilde{\omega}
 \\&  =:A_1+A_2+A_3 +A_4+A_5  .\nonumber
\end{align*}
\normalsize
We consider first the last three terms, the simplest ones.
By \eqref{eq:discrest2} we have
\begin{align*}
|A_3|&\lesssim
 \| {\eta}\chi_A\|_{L^1}|z | \   |\dot{z}  - \widetilde{z}| \lesssim     A ^{\frac{3}{2}}   A^{-1} \| \sech \( \frac{2}{A}x\) {\eta}\|_{L^2}  |\dot{z}  - \widetilde{z}|
 \\& \lesssim  A ^{\frac{1}{2}} \delta \(\| \eta \| _{\Sigma _A}  ^2+ \| \eta \| _{\widetilde{\Sigma}  }^{2} +|z|^4 \) .
\end{align*}
Since $\|  \partial _\omega g ^{(\omega)} \|_{L^\infty}\lesssim \<x \>$, using also \eqref{eq:discrest1}  we have
  \begin{align*}
|A_4|&\lesssim A
 \| {\eta}\chi_A\|_{L^1}|z | ^4   |\dot{\omega}  - \widetilde{\omega}| \lesssim     A ^{\frac{5}{2}}  |z | ^4  \| \eta \| _{\Sigma _A}   |\dot{\omega}  - \widetilde{\omega}|   \lesssim  {\delta}   A ^{\frac{5}{2}} |z | ^4
   .
\end{align*}
Finally, using \eqref{eq:estpar}  we have
\begin{align*}
|A_5|&\lesssim A
 \| {\eta}\chi_A\|_{L^1}|z | ^2   |  \widetilde{\omega}| \lesssim     A ^{\frac{5}{2}}  |z | ^4  \| \eta \| _{\Sigma _A}     \lesssim  {\delta}   A ^{\frac{5}{2}}  |z | ^4 .
\end{align*}
Turning to the main terms,
we have
\begin{align*}
   A_2 & =  \omega \<  J \eta , 2\im \lambda (p  ) \chi_A \( z^2  g ^{(\omega)}- \overline{z} ^2     \overline{g} ^{(\omega)}\)            \> + \<  J    \eta , \chi_A \( 2z     \widetilde{z}_{\mathcal{R}}  g ^{(\omega)}+  2\overline{z}  \ \overline{\widetilde{z}}_{\mathcal{R}} \ \overline{g} ^{(\omega)}\)            \> = A _{21}+ A _{22}.
\end{align*}
By \eqref{eq:estpar}   proceeding like for $A_3$,
\begin{align*}
|A _{22}|&\lesssim
 \| {\eta}\chi_A\|_{L^1}|z | ^3   \lesssim   { \delta }  A ^{\frac{3}{2}}    \| \eta \| _{\Sigma _A}  |z| ^2 \lesssim   \delta  A ^{\frac{3}{2}}
 \( \| \eta \| _{\Sigma _A}^{2}+  |z|^{4} \) .
\end{align*}
By \eqref{eq:nls4} and  by \eqref{eq:eqsatg2}
for the cancellation,  we have
  \small\begin{align*}&
A_1+A _{21}=  -\xcancel{
\<  J
	 {\eta},\chi_A   \mathcal{L}_{\omega}  (  {z}^{2} {g} ^{(\omega)} + \overline{{z}}^{2} \overline{g}^{(\omega)} )\> }  + \xcancel{A _{21}}  +  ( \widetilde{\vartheta}  _{\mathcal{R}} + \widetilde{\vartheta} - \dot \vartheta)
\< \eta ,  \chi_A   (  {z}^{2} {g} ^{(\omega)} + \overline{{z}}^{2} \overline{g}^{(\omega)} )\>
\\& + \dot D \<\partial _x \eta ,   \chi_A   (  {z}^{2} {g} ^{(\omega)} + \overline{{z}}^{2} \overline{g}^{(\omega)} )\>
-\< J   e^{J \vartheta +D\partial _x}D_\Theta \phi [\Theta]  (\dot \Theta  -\widetilde{\Theta})  , \chi_A   (  {z}^{2} {g}^{(\omega)}  + \overline{{z}}^{2} \overline{g} ^{(\omega)}) \>
\\
	& + \< (2\chi _A '\partial _x + \chi _A '')
	 {\eta} , {z}^{2} {g} ^{(\omega)} + \overline{{z}}^{2} \overline{g}^{(\omega)} \>  +\<  \(D f( \phi  [\omega ,v,  z  ] -D f( \phi  [\omega ,  v  ] )\)  \eta  ,\chi_A   (  {z}^{2} {g} ^{(\omega)} + \overline{{z}}^{2} \overline{g}^{(\omega)}  ) \>  \\& +\<  f( \phi [\omega ,v , z  ] + \eta  )  -  f( \phi  [\omega ,v , z  ]  ) -D f( \phi  [\omega ,  v  ]  ) \eta  ,\chi_A   (  {z}^{2} {g} ^{(\omega)} + \overline{{z}}^{2} \overline{g}^{(\omega)} ) \>
\\&+\< R [\omega ,v , z ],\chi_A(  {z}^{2} {g}^{(\omega)}  + \overline{{z}}^{2} \overline{g} ^{(\omega)})) \>
=A_{11}+A_{12}+A_{13}+A_{14} + A_{15}+A_{16} .
\end{align*}
\normalsize
It is easy to see the following
\begin{align*}
   \sum _{j=1}^{3}|A _{1j}| \lesssim  {\delta} A ^{\frac{3}{2}}
  \(  |z|^4   + \| \eta \| _{\Sigma _A}^{2}+ \| \eta \| ^{2}_{\widetilde{\Sigma}}  \)  .
\end{align*}
We have
\begin{align*}
   A_{14} &= \<  \(D f( \phi  [\omega ,v,  z  ] -D f(  \phi  [\omega ,v   ]  )\)  \eta  ,\chi_A   (  {z}^{2} {g} ^{(\omega)} + \overline{{z}}^{2} \overline{g} ^{(\omega)} ) \> \\& + \<  \(D f( \phi  [\omega ,v   ] -D f( \phi   _{\omega}  )\)  \eta  ,\chi_A   (  {z}^{2} {g} ^{(\omega)} + \overline{{z}}^{2} \overline{g}^{(\omega)}  ) \> =:  A_{141} + A_{142}\end{align*}
   On the support of $\chi _A$ we can assume $|\phi  [\omega ,v   ] | \gg |\widetilde{\phi}  [\omega ,v,  z  ] |$.
   We can write
  \begin{align*}
    A_{141} &= \int _0 ^1 \< D^2f(\phi  [\omega ,v   ] +s  \widetilde{\phi}  [\omega ,v,  z  ]  ) ( \widetilde{\phi}  [\omega ,v,  z  ]  ,\eta  )  ,\chi_A   (  {z}^{2} {g} ^{(\omega)} + \overline{{z}}^{2} \overline{g} ^{(\omega)} ) \>
  \end{align*}
   with
   \begin{align*}
      |  A_{141}|&\lesssim   |z|^2 \< \phi _\omega ^{p-2} |\widetilde{\phi}  [\omega ,v,  z  ]| \ |\eta| , \chi_A\> \lesssim
      |z|^3 \int _{\R}  e^{|x| \sqrt{\omega}  \( 2-p- \sqrt{1-\lambda (p)} \)  } |\eta | dx \lesssim |z|^3 \| \| \eta \| _{\widetilde{\Sigma}}.
   \end{align*}
   We have
    \begin{align*}
   A_{142}=
   \< \( e^{\frac{\im}{2} v x} D f( \phi   _{\omega}  )e^{-\frac{\im}{2} v x}- D f( \phi   _{\omega}  )\)  \eta  ,\chi_A   (  {z}^{2} {g} ^{(\omega)} + \overline{{z}}^{2} \overline{g} ) \> .
\end{align*}
So from this we get
\begin{align*}
  |  A_{142}|&\lesssim  |z|^2  |v|\<  |x| \phi _\omega ^{p-1},  |\eta | \>\lesssim  A ^{3/2}|v| |z|^2 \| \eta \| _{\widetilde{\Sigma}}\lesssim \delta  A ^{3/2}  \(  |z|^4    + \| \eta \| ^{2}_{\widetilde{\Sigma}}  \)  .
\end{align*}
Since we can assume $|{\phi}  [\omega ,v   ]|\gg |\widetilde{{\phi}}  [\omega ,v  ,z ]| +|\eta |$  in $\supp   \chi _A$  we write
\begin{align*}
   A_{15}&=\int  _{[0,1] ^2}  \<  D^2f( \phi [\omega ,v  ] +     \tau  ( \widetilde{\phi } [\omega ,v , z  ]+    s\eta )  )   \( \eta  , \widetilde{\phi } [\omega ,v , z  ]+    s\eta \) ,\chi_A   (  {z}^{2} {g} ^{(\omega)} + \overline{{z}}^{2} \overline{g}^{(\omega)} ) \> ds \ d\tau
\end{align*}
so that
\begin{align*}
  |A_{15}|&\lesssim |z|^2\int _{\R}\chi _A |\phi _\omega|^{p-2} |\eta |  \( |\widetilde{\phi } [\omega ,v , z  ]|  +|\eta |\) dx  \\& \lesssim  |z|^3 \int _\R    e^{|x| \sqrt{\omega}  \( 2-p- \sqrt{1-\lambda (p)} \)  } |\eta | dx + |z|^2e ^{2A (2-p)}A^2 \| \eta \|_{\Sigma _A}^2
\end{align*}
We have
\begin{align*}
  A_{16} &= \< R [\omega   , z ],\chi_A(  {z}^{2} {g}^{(\omega)}  + \overline{{z}}^{2} \overline{g} ^{(\omega)})) \> \\& +\<  \( e^{\frac{\im}{2} vx} -1\)   R [\omega  , z ],\chi_A(  {z}^{2} {g}^{(\omega)}  + \overline{{z}}^{2} \overline{g} ^{(\omega)})) \> =:  A_{161}+ A_{162}.
\end{align*}
Then we have
\begin{align*}
  | A_{162}|\lesssim |v| A |z|^2 \| R [\omega  , z ]\| _{L^1 \( \R \)} \lesssim |v| A |z|^4 \lesssim \delta A |z|^4.
\end{align*}
By \eqref{def:phiomegaz} and the line below \eqref{eq:eqsatg2}   we have
\begin{align}\label{eq:A161}
   A_{161}= \omega ^{\frac{p+1}{2(p-1)}}\<  R [  z ],\chi  _{A\sqrt{\omega}} (  {z}^{2} {g}   + \overline{{z}}^{2} \overline{g})  \> .
\end{align}
Omitting the harmless factor $\omega ^{\frac{p+1}{2(p-1)}}$ we have
\small
\begin{align*}
  A_{161} &=  \< f(\phi [  {z}]) -  f(\phi )          - Df(\phi   )\widetilde{\phi} [  {z}]   ,\chi  _{A\sqrt{\omega}}  (  {z}^{2} {g}  + \overline{{z}}^{2} \overline{g}  ) \> \\& +
 \<   - \widetilde{\vartheta}_\mathcal{R} \phi [  {z}]+ \im \widetilde{\omega}_\mathcal{R} \Lambda_p\phi [  {z}]   -\im \widetilde{D}_\mathcal{R} \partial _x \phi  [  {z}]-\frac{\widetilde{v}_\mathcal{R} x}{2} \phi  [  {z}]+\im  D_{z}\phi  [  {z}]\widetilde{z }_\mathcal{R}   ,\chi  _{A\sqrt{\omega}}  (  {z}^{2} {g}  + \overline{{z}}^{2} \overline{g}  )\> \\& =: A_{1611} + A_{1612}.
\end{align*}\normalsize
We claim
\begin{align}\label{eq:A152}
  | A_{1612}| \lesssim  A ^{-1}   |z|^4.
\end{align}
Indeed, for example we have
\begin{align*}
  &  \widetilde{\vartheta}_\mathcal{R}    \< \phi [  {z}],  \chi  _{A\sqrt{\omega}} (  {z}^{2} {g}  + \overline{{z}}^{2} \overline{g}  ) \> =  \widetilde{\vartheta}_\mathcal{R}    \<    \bcancel{ \phi  }    +   \widetilde{\phi} [  {z}]  ,     {z}^{2} {g}  + \overline{{z}}^{2} \overline{g}   \> \\& -  \widetilde{\vartheta}_\mathcal{R}     \<  \phi [  {z}], \(1 -  \chi  _{A\sqrt{\omega}}  \) (  {z}^{2} {g}  + \overline{{z}}^{2} \overline{g} )\>    = O\( z^5 \) + O\(  e ^{-\kappa \omega _0 A }z^4 \)
\end{align*}
where we used the orthogonality \eqref{eq:dirsum}--\eqref{eq:dirsum1}, the bound \eqref{eq:estpar} and the exponential decay of $\xi$. The other terms forming $A_{1612}$ can be bounded similarly, yielding
\eqref{eq:A152}.

  Exploiting   $\phi \gg |\widetilde{\phi} [  {z}]|$ on $\supp  \chi  _{A\sqrt{\omega}}$  we write  the Taylor formula   \small
\begin{align*}
  &A_{1611}  =  2^{-1} \<   D^2f(\phi )\(  \widetilde{\phi} [  {z}]\) ^{2}  ,\chi  _{A\sqrt{\omega}}(  {z}^{2} {g}   + \overline{{z}}^{2} \overline{g}  )) \> \\& +\int _{0}  ^{1} \frac{(1-s)^2}{2}
  \<     D^3f(\phi  +  s  \widetilde{\phi} [  {z}] )  \(   \widetilde{\phi} [  {z}] \) ^{3}  ,\chi _{A\sqrt{\omega}}(  {z}^{2} {g}   + \overline{{z}}^{2} \overline{g}  ) \> ds
  =:  A_{16111} + A_{16112}.
\end{align*}\normalsize
Then
\begin{align*}
  |A_{16112}|&\lesssim   |z|^2   \int _{|x|\le 2A\sqrt{\omega} }  |\phi | ^{p-3}  |\widetilde{\phi} [  {z}]|^{3} dx \lesssim  C _{A,\omega } |z|^5    \text{ where }  C _{A,\omega }=
  \int _{|x|\le 2A\sqrt{\omega} }   e^{|x| \( 3-p-3\sqrt{1-\lambda (p)}  \)}   dx.
\end{align*}
Here we take   $\delta>0$   so that
\begin{align*}
  |A_{16112}|\lesssim  C _{A,\omega } \delta  |z|^4 \le A ^{-1}   |z|^4.
\end{align*}
We write
\begin{align*}
   A_{16111} &=  2^{-1} \<   D^2f(\phi )\(  \widetilde{\phi} [  {z}]\) ^{2}  , {z}^{2} {g}   + \overline{{z}}^{2} \overline{g}    \> \\&  2^{-1} \<   D^2f(\phi )\(  \widetilde{\phi} [  {z}]\) ^{2}  , (1- \chi  _{A\sqrt{\omega}}) \( {z}^{2} {g}   + \overline{{z}}^{2} \overline{g} \)   \> =:  A_{161111}+ A_{161112}.
\end{align*}
Then
\begin{align*}
  | A_{161112}|\lesssim   |z|^4 \int _{|x|\ge \frac{A\sqrt{\omega}}{2}} e^{|x| \( 2-p- 2\sqrt{1-\lambda (p)}  \)}  dx  \lesssim A ^{-1} |z|^4.
\end{align*}
By    \eqref{eq:derf2}   for $\xi = (\xi _1, \xi _2) ^ \intercal$,  $X=\widetilde{\phi} [  {z}] =(z\xi _1 + \overline{z}\overline{\xi} _1)+\im (z\xi _2 + \overline{z}\overline{\xi} _2)$, $u= \phi  $ and identifying $\C=\R ^2$, we have
\begin{align} \label{eq:fgrhess}
  D^2f(\phi  )\(  \widetilde{\phi} [  {z}] \) ^{2} = (p-1)\phi   ^{p-2}  \begin{pmatrix} p (z\xi _1 + \overline{z}\overline{\xi} _1)^2+(z\xi _2 + \overline{z}\overline{\xi} _2) ^2
	    \\ 2(z\xi _1 + \overline{z}\overline{\xi} _1) (z\xi _2 + \overline{z}\overline{\xi} _2)
\end{pmatrix} .
\end{align}
 Then
 we have
\begin{align*}
   A_{161111} & =    (p-1) |z|^4  \gamma (p )    \\& +  2(p-1) |z|^2 \( \<    \phi _\omega ^{p-2}  \( p |\xi _1|^2+ |\xi _2|^2 \) , z^2 g_1 ^{(\omega)}\> + 2\<    \phi _\omega ^{p-2}    (\xi _1  \overline{\xi} _2  + \overline{\xi} _1   {\xi} _2)   ,z^2  g_2^{(\omega)} \> \) \\& +   (p-1) \( \<    \phi _\omega ^{p-2}  \( p  \overline{\xi} _1 ^2+  \overline{\xi} _2 ^2 \) , z^4 g_1^{(\omega)} \> +2 \<    \phi _\omega ^{p-2}     \overline{\xi} _1  \overline{\xi} _2     ,z^4  g_2^{(\omega)} \> \)   \\&  =:  A_{161111}+ A_{161112} +A_{161113} .
\end{align*}
We claim that We claim that  there exists a function $\mathcal{I}_{\mathrm{FGR}}$ like in the statement such that
\begin{align}\label{eq:classnormform}
 \left |    - \frac{d}{dt} \mathcal{I}_{\mathrm{FGR}}  +  \sum _{j=2}^{3}  A_{16111j}     \right |\lesssim |z|^5+ |z | ^3   |\dot \Theta - \widetilde{\Theta}|.
\end{align}
To   $n=2,4$   write
\begin{align}\label{eq:nform}
   \frac{d}{dt}z  ^n&= n   \widetilde{ {z}}  z   ^{n-1} +n  \( \dot z -  \widetilde{ {z}}   \) z   ^{n-1}   =   n\im \lambda        z ^n +    n   \widetilde{ {z}}_{\mathcal{R}}  z   ^{n-1}     +n  \( \dot z -  \widetilde{ {z}}   \) z   ^{n-1}.
\end{align}
Then for example
\begin{align*}&
  \<    \phi   ^{p-2}     \overline{\xi} _1  \overline{\xi} _2     ,z^4  g_2  \> = \<    \phi   ^{p-2}     \overline{\xi} _1  \overline{\xi} _2     ,   \frac{1}{4\im \lambda (p) }
    \(  \frac{d}{dt} z^4  - 4  \widetilde{ {z}}_{\mathcal{R}}  z   ^{3} -4  \( \dot z -  \widetilde{ {z}}   \) z   ^{3} \)  g_2  \>
\end{align*}
and applying the Leibnitz rule for the time derivative, it is easy to obtain the claim.
Summing up  and recalling also that we have to  reintroduce the      factor $\omega ^{\frac{p+1}{2(p-1)}}$  of \eqref{eq:A161} which we have omitted for simplicity,      we obtain   \eqref{eq:lem:FGR11}.

\qed

\textit{Proof of Proposition \ref{prop:FGR}.
}
Integrating \eqref{eq:lem:FGR11}  we obtain $\| z^2 \| _{L^2(I)}   ^2 \lesssim \sqrt{\delta} + A ^{-1}\epsilon ^2$ yielding \eqref{eq:FGRint}.

\qed

\section{High energies: proof of Proposition \ref{prop:1virial}}\label{sec:vir}

Following  the framework in Kowalczyk et al. \cite{KMM3}  and
using the function $\chi$  in \eqref{eq:chi}
we  consider the   function
\begin{align}\label{def:zetaC}
\zeta_A(x):=\exp\(-\frac{|x|}{A}(1-\chi(x))\)   \text {  and   }  \varphi_A (x):=\int_0^x \zeta_A^2(y)\,dy
\end{align}
  and  the vector field
\begin{align}\label{def:vfSA} \  S_A:=   \varphi_A'+2 \varphi   _{ A}\partial_x .
\end{align}
Next, following Martel  \cite{Martel1} we introduce   \begin{align}\label{def:funct:vir}
\mathcal{I} := 2^{-1}\<  \im \xi   ,S_A  \xi  \>     \text{ where }  \xi :=e^{-\frac{\im}{2} xv}  \eta .
\end{align}
\begin{lemma}\label{lem:1stV1}
There exists a fixed constant $C>0$ s.t. for an arbitrary small number
\begin{align} &
\|   \eta  \|_{\Sigma _A}^2   \le C\left [  \dot{\mathcal{I}} +  \|    \eta \| _{\widetilde{\Sigma}} ^2  +|\dot \Theta   - \widetilde{\Theta}| ^2+     |z |^4 \right ].   \label{eq:lem:1stV1}
\end{align}
\end{lemma}
\proof We have
\begin{align}\nonumber
  \dot {\mathcal{I}} &= - \<  \dot \xi ,    \im S_A\xi  \>   = -\<  \dot \eta -\frac{\im }{2} \dot v x \eta , \im S_A \eta +v \varphi   _{ A} \eta \>
  \\&=-\<  \dot \eta,    \im S_A\eta \> -\frac{\dot v }{4} \< [S_A, x]\eta , \eta \> -v\<  \dot \eta , \varphi   _{ A} \eta \> .\label{eq:vir2}
\end{align}
From \eqref{eq:nls3}   we obtain
\begin{align*}
   -\<  \dot \eta,    \im S_A\eta \>  &= -\<  \partial _x^2    \eta ,S_A\eta  \> -  \<    R [\omega ,v , z ] ,S_A\eta  \>   - \frac{1}{2}   \dot D   \< \eta  ,\im [\partial _x ,S_A]\eta  \>  \\& -  \<    f( \phi  [\omega ,v , z ]+ \eta  )  -  f( \phi  [\omega ,v , z ]   )     ,S_A\eta  \>   + O\(  |\dot \Theta   - \widetilde{\Theta}|   \|    \eta \| _{\widetilde{\Sigma}}  \)
\end{align*}
and
\begin{align*}
   -v\<  \dot \eta , \varphi   _{ A} \eta \>  &= -v\< \im  \partial _x^2    \eta , \varphi   _{ A}  \eta  \>  -v\dot D \< \eta ', \varphi _A \eta \> -v\< \im  R[\omega ,v   , z ] , \varphi   _{ A} \eta  \>
   \\& -v\<  \im  \( f( \phi [\omega ,v  , z ] + \eta  )  -  f( \phi  [\omega ,v, z ]  )  \) , \varphi   _{ A} \eta  \> .
\end{align*}
From  \cite{Martel1} we have
\begin{align}
  - \<  \partial _x^2    \eta ,S_A\eta  \> \ge   2 \|   (\zeta _A \eta)' \| _{L^2}^2- \frac{C }{A}  \|    \eta \| _{\widetilde{\Sigma}} ^2.\label{eq:poscomm}
\end{align}
We have
\begin{align*}
  -v\dot D \< \eta ', \varphi _A \eta \> =  \frac{v}{2}\dot D  \| \zeta  _A \eta  \| _{L^2} ^2 =   \frac{v^2}{2}   \| \zeta  _A \eta  \| _{L^2} ^2 + \frac{v}{2}(\dot D -v)  \| \zeta  _A \eta  \| _{L^2} ^2\ge   \frac{v}{2}(\dot D -v)  \| \zeta  _A \eta  \| _{L^2} ^2.
\end{align*}
So we obtain
\begin{align}\label{eq:vir3}
   2 \|   (\zeta _A \eta)' \| _{L^2}^2  &\le  \dot {\mathcal{I}}  + \frac{1}{2}   \dot D   \< \eta  ,\im [\partial _x ,S_A]\eta  \>  +v\< \im  \partial _x^2    \eta , \varphi   _{ A}  \eta  \>
   \\& \label{eq:vir4} + \<    R [\omega ,v , z ] ,S_A\eta  \>  +v\< \im  R[\omega ,v   , z ] , \varphi   _{ A} \eta  \>
   \\& \label{eq:vir6} +   \<    f( \phi  [\omega ,v , z ]+ \eta  )  -  f( \phi  [\omega ,v , z ]   )     ,S_A\eta  \>       \\& +v\<  \im  \( f( \phi [\omega ,v  , z ] + \eta  )  -  f( \phi  [\omega ,v, z ]  )  \) , \varphi   _{ A} \eta  \>  \label{eq:vir70}\\&\label{eq:vir7}+\frac{\dot v }{4} \< [S_A, x]\eta , \eta \> -\frac{v}{2}(\dot D -v)  \| \zeta  _A \eta  \| _{L^2} ^2 + O\(  |\dot \Theta   - \widetilde{\Theta}|   \|    \eta \| _{\widetilde{\Sigma}}  \)
\end{align}
 We  have
\begin{align*}
     \frac{1}{2}   \dot D   \< \eta  ,\im [\partial _x ,S_A]\eta  \> &=   \dot D   \< \eta  ,\im [\partial _x ,\varphi _A]\eta'  \> =  \dot D   \< \eta  ,\im \zeta _A^2    \eta'  \>
\end{align*}
and
\begin{align*}
   v\< \im  \partial _x^2    \eta , \varphi   _{ A}  \eta  \> =-v\< \im      \eta  ',\zeta _A ^2  \eta  \>
\end{align*}
and so  by \eqref{eq:oertstab1} and     \eqref{eq:discrest1}
\begin{align}\label{eq:vir8}
  \left |    \frac{1}{2}   \dot D   \< \eta  ,\im [\partial _x ,S_A]\eta  \> +v\< \im  \partial _x^2    \eta , \varphi   _{ A}  \eta  \> \right |=  \left | (\dot D -v)  \< \eta  ,\im \zeta _A^2    \eta'  \> \right | \lesssim  \delta ^2 \(   |z|^4+    \|   \eta \| _{\widetilde{\Sigma}}     ^2  +\| \eta \| _{\Sigma_{A}}^2\).
\end{align}
By \eqref{estR}, for a preassigned small constant $\delta _2>0$
\begin{align}
  | \<    R  [\omega ,v , z ],S_A\eta  \>|  +|v\< \im  R[\omega ,v   , z ] , \varphi   _{ A} \eta  \>|   \lesssim |z|^4 +\| \eta \| _{\widetilde{\Sigma}} ^2 +\delta _2 \| \eta \| _{ {\Sigma}_A} ^2.\label{eq:vir9}
\end{align}
By  \eqref{eq:oertstab1} and     \eqref{eq:discrest1}  we have
\begin{align}
  \left | \frac{\dot v }{4} \< [S_A, x]\eta , \eta \>  \right |&\lesssim   A |\dot v|   \| \eta \| _{L^2} ^2      \lesssim  A\delta ^2     \(   |z|^4+    \|   \eta \| _{\widetilde{\Sigma}}     ^2  +\| \eta \| _{\Sigma_{A}}^2\) \label{eq:vir10}
\end{align}
and
\begin{align}
  \left | \frac{v}{2}(\dot D -v)  \| \zeta  _A \eta  \| _{L^2} ^2  \right |&\lesssim \delta ^4         \(   |z|^4+    \|   \eta \| _{\widetilde{\Sigma}}     ^2  +\| \eta \| _{\Sigma_{A}}^2\) . \label{eq:vir11}
\end{align}
We write \begin{equation*}
      \<   f( \phi  [\omega  ,v, z ]+ \eta  )  -  f( \phi  [\omega ,v , z ]   )   - f\left(  \eta  \right) ,S_A {\eta} \> + \<  f\left(  \eta  \right) ,S_A {\eta} \>=:B_{ 1}+B_{ 2}.
 \end{equation*}
Then
\begin{align*}
  B_{ 2} =  \<  f\left(  \eta  \right)   \overline{{\eta}} - F\left(  \eta  \right) , \zeta _A ^2 \>  =   \frac{p}{p+1} \int _\R |\eta | ^{p+1} \zeta _A ^2 dx
\end{align*}
and we use   the crucial estimate by    Kowalczyk et al. \cite{KMM3}
\begin{align}\label{eqpurepower}
  |B_{ 2}| \lesssim   A^2 \| \eta \| ^{p-1}_{L^\infty (\R )}   \| (\zeta _A\eta ) ' \| ^{2}_{L^2 (\R )} \lesssim A ^{-1} \| (\zeta _A\eta ) ' \| ^{2}_{L^2 (\R )}.
\end{align}
\begin{claim}
  \label{claim:estB1} We have  for any small preassigned $\delta _2>0$
  \begin{align}\label{eq:claim:estB1}
  |B_{ 1}| \lesssim    \|    \eta \| _{\widetilde{\Sigma}} ^2 + \delta _2  \|   (\zeta _A \eta)' \| _{L^2}^2 .
\end{align}
\end{claim}
\proof We consider
\begin{align}& \label{eq:claim:estB11}
   \Omega_{1,t  }=\{x\in \R\ |\  |   \eta(t,x)|\leq 2|\phi[\omega (t) ,v(t),   z(t)]|\}  \text{ and }\\&  \Omega_{2,t  }=\R\setminus \Omega_{1,t  }=\{x\in \R\ |\ |  \eta(t,x)|> 2|\phi[\omega (t) ,v(t),  z(t)]|\} , \nonumber
\end{align}
and split accordingly
\begin{align*}
  B_1 =\sum _{j=1,2} \<   f( \phi  [\omega ,v , z ]+ \eta  )  -  f( \phi  [\omega ,v , z ]   )   - f\left(  \eta  \right) ,S_A {\eta} \> _{L^2( \Omega_{j,t  })} =:B _{11}+B _{12}.
\end{align*}
We have
\begin{align*}
  |B _{11}|&\le \int _0 ^1 \int _{\Omega_{1,t  }} \left | \<   \( D   f(  \phi  [\omega ,v , z ]+ s\eta  )  -  Df(s \eta   )\)  \eta  ,S_A {\eta}           \> _\C  \right | dx ds\\& \lesssim
  \int _{\Omega_{1,t  }}   \zeta _A ^2|\phi  [\omega  , z ]  | ^{p -1} |\eta | ^2 dx +  \int _{\Omega_{1,t  }} \zeta _A ^{-1}  |\varphi  _A |  |\phi  [\omega  , z ]  | ^{p -1}  |\eta |    (\zeta _A \eta) ' - \zeta _A '\eta   |      dx =: B _{111}+B _{112}.
\end{align*}
Then
\begin{align*}
  B _{111} \lesssim  \|    \eta \| _{\widetilde{\Sigma}} ^2
\end{align*}
and, for a fixed small and preassigned $\delta _2>0$,
\begin{align*}
  B _{112} \lesssim  \int _{\Omega_{1,t  }}  \zeta _A ^{-1}  |\psi _A |  |\phi  [\omega  , z ]  | ^{p -1}  |\eta |   \(   |(\zeta _A \eta) '|  + \frac{1}{A} \zeta _A |\eta |  \)         dx  \lesssim  \|    \eta \| _{\widetilde{\Sigma}} ^2 + \delta _2  \|   (\zeta _A \eta)' \| _{L^2}^2 .
\end{align*}
We have
\begin{align*}
  |B _{12}|&\le \int _0 ^1 \int _{\Omega_{2,t  }} \left | \<   \( D   f( s \phi  [\omega ,v , z ]+ \eta  )  -  Df(s \phi  [\omega ,v , z ]   )\)  \phi  [\omega ,v , z ]  ,S_A {\eta}           \> _\C  \right | dx ds\\& \lesssim
  \int _{\Omega_{2,t  }}   \zeta _A ^2|\eta | ^{p } |\phi  [\omega  , z ]| dx +  \int _{\Omega_{2,t  }} \zeta _A ^{-1}  |\varphi _A |  |\eta | ^{p-1 } |\phi  [\omega  , z ]|  | (\zeta _A \eta) ' - \zeta _A '\eta   |      dx =: B _{121}+B _{122}.
\end{align*}
We have
\begin{align*}
  B _{121} \lesssim  \int _{\Omega_{2,t  }}   \zeta _A ^2|\eta | ^{p } |\phi  [\omega  , z ]| ^{2-p+p-1}    dx \le  \int _{\R}   |\eta | ^{2 } |\phi  [\omega  , z ]| ^{p-1}    dx \lesssim \|    \eta \| _{\widetilde{\Sigma}} ^2.
\end{align*}
We have similarly, for a fixed small and preassigned $\delta _2>0$ and completing the proof of Claim \ref{claim:estB1},
\begin{align*}
  B _{122} \lesssim  \int _{\Omega_{2,t  }}  \zeta _A ^{-1}  |\psi _A | |\eta |   |\phi  [\omega  , z ]| ^{p-1}  \(   |(\zeta _A \eta) '|  + \frac{1}{A} \zeta _A |\eta |  \)         dx  \lesssim  \|    \eta \| _{\widetilde{\Sigma}} ^2 + \delta _2  \|   (\zeta _A \eta)' \| _{L^2}^2 .
\end{align*}
\qed

We consider the term in \eqref{eq:vir70} which equals
\begin{align*}
 \sum _{j=1,2}  v\<  \im  \( f( \phi [\omega ,v  , z ] + \eta  )  -  f( \phi  [\omega ,v, z ]  ) -f(   \eta  )  \) , \varphi   _{ A} \eta  \> _{L^2( \Omega_{j,t  })} =:C_{1 }+C _{2 }
\end{align*}
and then we proceed like in  Claim \ref{claim:estB1}.
We have
\begin{align*}
  |C _{1 }|&\le |v|\int _0 ^1 \int _{\Omega_{1,t  }} \left | \<  \im \( D   f(  \phi  [\omega ,v , z ]+ s\eta  )  -  Df(s \eta   )\)  \eta  ,\varphi   _{ A} {\eta}           \> _\C  \right | dx ds\\& \lesssim
 |v| \int _{\Omega_{1,t  }}   |\varphi _A |  \ |\phi  [\omega  , z ]  | ^{p -1} |\eta | ^2 dx  \lesssim  \delta \|    \eta \| _{\widetilde{\Sigma}} ^2 .
\end{align*}
We have
\begin{align*}
  |C _{ 2}|&\le |v| \int _0 ^1 \int _{\Omega_{2,t  }} \left | \<   \( D   f( s \phi  [\omega ,v , z ]+ \eta  )  -  Df(s \phi  [\omega ,v , z ]   )\)  \phi  [\omega ,v , z ]  ,\varphi   _{ A}  {\eta}           \> _\C  \right | dx ds\\& \lesssim
 |v| \int _{\Omega_{2,t  }}     |\varphi _A |  \     |\eta | ^{p } |\phi  [\omega  , z ]| ^{2-p+p-1}    dx \le |v| \int _{\R}   |\eta | ^{2 } |\phi  [\omega  , z ]| ^{p-1}    dx \lesssim  \delta   \|    \eta \| _{\widetilde{\Sigma}} ^2.
\end{align*}
 Collecting together \eqref{eq:vir3}--\eqref{eq:vir7}
 and the estimates we get for a small constant $\delta _2>0$
  \begin{align*}
     2 \|   (\zeta _A \eta)' \| _{L^2}^2 \le   \dot {\mathcal{I}}  +\delta _2 \|   (\zeta _A \eta)' \| _{L^2}^2 + O\( |z|^4+   \|    \eta \| _{\widetilde{\Sigma}} ^2 \)
  \end{align*}
 which yields a bound of $\|   (\zeta _A \eta)' \| _{L^2}^2$. Finally to get  \eqref{eq:lem:1stV1} we use the following,  see Lemma 6.2 \cite{CMS23},  \begin{align*} \| \sech  \(   \frac{2}{A}x\) \eta '\|_{L^2}^2 + A^{-2}\| \sech  \(   \frac{2}{A}x\) \eta \|_{L^2}^2 \lesssim   \|   ( \zeta _A \eta )' \| _{L^2(\R )} ^2    +A ^{-1}  \|    \eta \| _{\widetilde{\Sigma}}  ^2.\end{align*}
\qed

 \textit{Proof of Proposition \ref{prop:1virial}.}
Integrating  inequality \eqref{lem:1stV1} we obtain \eqref{eq:sec:1virial1}.

\qed

\section{Low energies: proof of Proposition \ref{prop:smooth11}}\label{sec:smooth1}

Applying $U^{-1}$ to equation \eqref{eq:nls4} we get
\small
\begin{align}\label{eq:nls5}
   \partial _t(  U^{-1} \eta )   &=  \im  H _{\omega} U^{-1} \eta  -  \im \sigma _3   ( \widetilde{\vartheta}  _{\mathcal{R}} + \widetilde{\vartheta} - \dot \vartheta) U^{-1}  \eta  +   \dot D    \partial _xU^{-1} \eta\\& \nonumber  -   e^{\im \sigma _3 \vartheta +D\partial _x} U^{-1}D_\Theta \phi [\Theta]  (\dot \Theta  -\widetilde{\Theta})   +\im \sigma _3   U^{-1} \( D f( \phi  [\omega ,v , z ]  ) - D f( \phi   _{ \omega }   )\) \eta
    \\&-\im \sigma _3   U^{-1} \( f( \phi [\omega ,v , z ] + \eta  )  -  f( \phi  [\omega ,v , z ]  ) -D f( \phi  [\omega ,v , z ]  ) \eta\) -\im \sigma _3   U^{-1}R[\omega ,v , z ].  \nonumber
\end{align}\normalsize
Set     $v:= \chi _B   U^{-1}\eta $  for $B$ like in \eqref{eq:relABg}.   We denote  $P_d(\omega) $ the discrete spectrum projection and   $P_c(\omega)$ the continuous spectrum projection associated to $H _{\omega}$, which are closely related to the corresponding projections for $\mathcal{L}_\omega$,   indeed  $ U P_a(H _{\omega})  U ^{-1} =  P_a(\mathcal{L}_{\omega})   $ for $a=c,d$.  Then we write
\begin{align}\label{eq:scomv} &   v=P_c(\omega _0) v + P_d(\omega _0) v   \text{ where we have}   \\&  P_d(\omega _0) v=  P_d(\mathcal{L}_{\omega}) \eta + U ^{-1} \( P_d(\mathcal{L}_{\omega _0}) - P_d(\mathcal{L}_{\omega}) \) \eta - P_d(\omega _0) \( 1 - \chi _B\)   U ^{-1} \eta \nonumber
\end{align}
     It is easy to check that
 \begin{align}
   \nonumber  \|  P_d(\omega _0) \( 1 - \chi _B\)   U ^{-1} \eta  \| _{L^{2,s}(\R ) } +    \|     U ^{-1} \( P_d(\mathcal{L}_{\omega _0}) - P_d(\mathcal{L}_{\omega}) \) \eta  \| _{L^{2,s}(\R ) } \le o_{B^{-1}}(1)   \| \eta \| _{\widetilde{\Sigma}}    \text{   for any $s\in \R$} .
 \end{align}
 Finally, with an argument similar to one in \cite[\S 7]{CM24D1}  by means of the orthogonality  in \eqref{61} it is easy to see that 
 \begin{align*}&
     \|     U ^{-1} P_d(\mathcal{L}_{\omega}) \eta   \| _{L^{2,s}(\R ) } \lesssim     |z|   \| \eta \| _{\widetilde{\Sigma}}    \text{   for any $s\in \R$} .
 \end{align*}
 So we conclude that  we have
\begin{align}
  \label{eq:scomv1}   \| P_d(\omega _0) v   \| _{L^{2,s}(\R ) } = o_{B^{-1}}(1)   \| \eta \| _{\widetilde{\Sigma}}    \text{   for any $s\in \R$} .
\end{align}
Setting $w= P_c(\omega _0)v$, we have      \small
\begin{align}\label{eq:nls6}
     &\partial _t w -  \im  H _{\omega _0} w=   -\im  \varpi  P_c(\omega _0)\sigma _3        w   +\im  P_c(\omega _0) \sigma _3 \( 2 \chi ' _B \partial _x + \chi _B '' \)  U^{-1}\eta  +     \dot D    P_c(\omega _0) \chi   _B \partial _xU^{-1} \eta \\& \label{eq:nls7}  -\im  \varpi  P_c(\omega _0)\sigma _3 P_d(\omega _0) v  + \im P_c (V _{\omega _{0}}    - V _{\omega  })w      + \im P_c (V _{\omega _{0}}    - V _{\omega  })    P_d(\omega _0) v   \\& \label{eq:nls8}  -  P_c(\omega _0) \chi _B e^{\im \sigma _3 \vartheta +D\partial _x} U^{-1}D_\Theta \phi [\Theta]  (\dot \Theta  -\widetilde{\Theta})   +\im P_c(\omega _0) \sigma _3 \chi _B   U^{-1} \( D f( \phi  [\omega ,v , z ]  ) - D f( \phi   _{ \omega }   )\) \eta
    \\&-\im P_c(\omega _0)  \chi _B\sigma _3   U^{-1} \( f( \phi [\omega ,v , z ] + \eta  )  -  f( \phi  [\omega ,v , z ]  ) -D f( \phi  [\omega ,v , z ]  ) \eta\) -\im P_c(\omega _0)\sigma _3  \chi _B  U^{-1}R[\omega ,v , z ]   \label{eq:nls9}
\end{align}\normalsize
where
\begin{align} \label{eq:nls61}
   \varpi := \widetilde{\vartheta}  _{\mathcal{R}} + \widetilde{\vartheta} - \dot \vartheta +\omega -\omega _0 .
\end{align}
Now we have
\begin{align}\nonumber
  w&=e^{\im t H _{\omega _0}}  w(0) - \im  \int _{0}^t e^{-\im (t-s) H _{\omega _0}}\varpi  (s) P_c(\omega _0)\sigma _3        w (s) ds \\& +  \im  \int _{0}^t e^{-\im (t-s) H _{\omega _0}}  P_c(\omega _0)\sigma _3         \( 2 \chi ' _B \partial _x + \chi _B '' \)  U^{-1}\eta (s) ds \nonumber \\& \nonumber  +    \int _{0}^t e^{-\im (t-s) H _{\omega _0}}   \dot D (s)   P_c(\omega _0) \chi   _B \partial _xU^{-1} \eta  (s) ds\\& +  \int _{0}^t e^{-\im (t-s) H _{\omega _0}} \text{lines \eqref{eq:nls7}--\eqref{eq:nls9}} ds.\label{eq:expv2}
\end{align}

The following is a simplified version of \cite[Lemma 7.2]{CM24D1}.

\begin{lemma}\label{lem:estw}
 For $S > 3/2$  we have
 \begin{equation}\label{eq:estw1}
   \| w \|  _{L^2(I , L ^{2,-S}(\R ) )}\le   o_{B^{-1}}(1)   \epsilon.
 \end{equation}

\end{lemma}  \proof
By  the analogue for $H_\omega $ of \eqref{eq:smooth111}, see \S 9 \cite{CM24D1},     we have
 \begin{align}
   \|    e ^{\im t H _{\omega _0}}   w(0) \| _{L^2(\R, L ^{2,-S}(\R ))} \lesssim   \|
   w(0)  \| _{L^2 (\R )} \lesssim \|
   \eta(0)  \| _{L^2 (\R )} \lesssim \delta \( = o_{B^{-1}}(1)   \epsilon  \).\nonumber
 \end{align}
Using the version for $H_\omega $ of Proposition \ref{prop:KrSch} (notice incidentally that in Krieger and Schlag  \cite{KrSch} the operators  are written as $H_\omega$)
we have
\begin{align*}&
   \|  \int _{0}^t e^{-\im (t-s) H _{\omega _0}}\varpi  (s) P_c(\omega _0)\sigma _3        w (s) ds \| _{L^2(I , L ^{2,-S}(\R ) )}  \lesssim  \| \varpi \| _{L^\infty(I)}   \|    w   \| _{L^2(I , L ^{2, S}(\R ) )} \\& \lesssim  \delta   \(\|    \chi _B   U^{-1}\eta   \| _{L^2(I , L ^{2, S}(\R ) )}     +  \|   [P_d(\omega _0), \chi _B]  U^{-1}\eta   \| _{L^2(I , L ^{2, S}(\R ) )}    \)  \\& \lesssim  \delta   \(   B ^s A A ^{-1}\|  \sech \( \frac{2}{A} x \)   \eta   \| _{L^2(I , L ^{2 }(\R ) )}     + o_{B^{-1}}(1)  \|   [ \eta   \| _{L^2(I , \widetilde{\Sigma} )}    \) \\& = o_{B^{-1}}(1) \( \|     \eta   \| _{L^2(I ,  \Sigma _A)} +    \|   [ \eta   \| _{L^2(I , \widetilde{\Sigma} )}  \) =  o_{B^{-1}}(1)   \epsilon.
\end{align*}
Using the analogue for $H_\omega $ of  \eqref{eq:smoothest11}  we consider
\begin{align*}&
   \|  \int _{0}^t e^{-\im (t-s) H _{\omega _0}}P_c(\omega _0)\sigma _3         \( 2 \chi ' _B \partial _x + \chi _B '' \)  U^{-1}\eta (s) ds    \| _{L^2(I , L ^{2,-S}(\R ) )} \\& \lesssim
  \|     \( 2 \chi ' _B \partial _x + \chi _B '' \)  \eta     \| _{L^2(I , L ^{2,\tau }(\R ) )}    \lesssim B ^{\tau-1} \| \sech \( \frac{2}{A}x\) \eta ' \| _{L^2(I ,L^2(\R ))} \\& + B ^{\tau -2} \| 1_{B\le |x|\le 2 B} \sech \( \frac{2}{A} x\)  \eta   \| _{L^2(I ,L^2(\R ))} \lesssim B ^{\tau-1} \|   \eta   \| _{L^2(I ,\Sigma _A )} \\& +  B ^{\tau -1} \( \left \|   \( \sech \( \frac{2}{A}x\) \eta  \) '    \right \| _{L^2(I ,L^2(\R ))}   +  \|    \eta \| _{L^2(I , \widetilde{\Sigma})}
  \) =o_{B^{-1}} (1)   \epsilon,
\end{align*}
where we used   $\tau \in (1/2,1) $   and, see   Merle and Raphael \cite[Appendix C]{MR4},
\begin{align*}
  \|   u \| _{L^2(|x|\le 2 B)}\lesssim  B \( \left \|  u '    \right \| _{ L^2(\R  )}   +  \|    u \| _{  \widetilde{\Sigma} }\)  .
\end{align*}
Next, using again   the version for $H_\omega $ of Proposition \ref{prop:KrSch}
\begin{align*}&
   \|   \int _{0}^t e^{-\im (t-s) H _{\omega _0}}   \dot D (s)   P_c(\omega _0) \chi   _B \partial _xU^{-1} \eta  (s) ds \| _{L^2(I , L ^{2,-S}(\R ) )}  \lesssim \| \dot D\|_{L^\infty (I)}\|     \chi   _B \partial _x  \eta    \| _{L^2(I , L ^{2, S}(\R ) )}  \\&\lesssim \delta B^S \| \sech \( \frac{2}{A} x\)  \eta '  \| _{L^2(I , L ^{2 }(\R ) )} \le o_{B^{-1}} (1) \| \eta \|  _{L^2(I ,\Sigma _A )} = o_{B^{-1}} (1)   \epsilon.
\end{align*}
Finally the terms in line \eqref{eq:expv2} can be similarly bounded using in particular the  analogue for $H_\omega $  of {Proposition} \ref{prop:KrSch}. The estimates are elementary and similar to \cite[Sect. 8]{CM2109.08108}.

\qed

\textit{Proof of Proposition \ref{prop:smooth11}.}
From \eqref{eq:scomv},  \eqref{eq:scomv1} and \eqref{eq:estw1} we have $\| v\| _{L^2(I , \widetilde{\Sigma})}\lesssim o_{B^{-1}}(1) \epsilon$.  Next, from  $v = \chi _B   U^{-1}\eta $,  and thanks to the relation  $A\sim B^3$  set in \eqref{eq:relABg}    we have
\begin{align} \nonumber
   \| \eta \| _{\widetilde{\Sigma}}&\lesssim  \| v \| _{\widetilde{\Sigma}} + \| (1- \chi _B   )\eta \| _{\widetilde{\Sigma}}\lesssim \| v \| _{\widetilde{\Sigma}} + A ^{-2} \|  \sech \( \frac{2}{A}x  \) \eta \| _{L^2}\\& \lesssim \| v \| _{\widetilde{\Sigma}} + A ^{-1} \|   \eta \| _{\Sigma _A}\label{eq:relAB0}
\end{align}
So by \eqref{eq:sec:1virial1}  and \eqref{eq:FGRint}  we get  the following, which implies \eqref{eq:sec:smooth11},
\begin{align} \nonumber
   \| \eta \| _{L^2(I,\widetilde{\Sigma})}&\lesssim    \| v \| _{L^2(I,\widetilde{\Sigma})} + A ^{-1} \|   \eta \| _{L^2(I, \Sigma _A)}\\& \lesssim  o_{B^{-1}} (1) \epsilon +  A ^{-1} \( \| z^2\|_{L^2(I)} + \| \eta \| _{L^2(I, \widetilde{\Sigma}   )}\) \lesssim   o_{B^{-1}}(1) \epsilon +  A ^{-1}\| \eta \| _{L^2(I, \widetilde{\Sigma}   )}    .  \nonumber
\end{align}

 \qed

\section{Proof of Proposition \ref{prop:FGR}}\label{sec:prop:FGR}

\textit{ Proof of Proposition \ref{prop:FGR}}  Here
    from  \eqref{eq:opH}    \begin{align*}
   &   H    =   \sigma _3 \(-\partial _x^2 +1\) +V (p,x)   , \\&
    V (p,x)   :=    M_0 \sech ^2   \(   \frac{p-1}2 x\)  \text{ with } M_0= -\( \frac{p+1}{2}
    \sigma _3      +\im  \frac{p-1}{2}  \sigma _2 \)   \frac {p+1}2   .  \nonumber
\end{align*}
From this   we see that $H$ is analytic of type $(A)$  in $p$  and this implies that $\lambda (p)$ depends analytically on $p \in (p_0, 3)$, see \cite[Ch. 12]{reedsimon}. By Lemma \ref{lem:simandmult}
 we already know that $\lambda (p)$ has (algebraic and geometric) multiplicity 1  as an eigenvalue of $H$  for $p\in (p_0, 3)$.  From analyticity in the sense of Kato, we conclude by  \cite[Theorem XII.8]{reedsimon} that  it is possible to define   $\xi$ dependent analytically in $p$. Proceeding like in \S 7  \cite{CM243}     it is possible to define the $g$   in \eqref{eq:eqsatg2}
to  dependent analytically in $p$. It follows that the function $\gamma (p) $ in \eqref{eq:fgrgamma} is analytic in  $p \in (p_0, 3)$. Finally, since it is possible to arrange things so that for
 $0<|p-3|\ll 1 $ we have $\gamma (p)\neq 0 $ as shown in \cite{CM24D1,CM242}, we conclude that  $(p_0, 3) \setminus \{   p\in (p_0, 3) : \quad \gamma(p)=0  \}$ is a discrete subset of $(p_0, 3)$, proving Proposition \ref{prop:FGR}.

  \qed

Finally we consider the following lemma, which was used in the proof of Lemma \ref{lem:simandmult} and where all the vector spaces are considered with corresponding scalar field $\C$.

\begin{lemma}\label{eq:multeig}
For any $p>1$ we have
  We have  $  \ker \(  \mathcal{L} -\im \lambda   \) = N_g\( \mathcal{L} -\im \lambda  \)$  if $\lambda   \in (0,1)$.
\end{lemma}
\proof Recall that $\mathcal{L} =J K$   where    $K= \diag \( L_+, L_- \)$.  The operator $K$ has a single and simple  negative eigenvalue and $\dim \ker  K =2$, see Weinstein \cite{W2}.
 Consider now
the Hermitian bilinear form $B(u)= \langle K
u , u \rangle $ using the Hermitian form $\< \cdot , \cdot \>$ in $L^2(\R , \C )$ introduced in \S \ref{sec:lin1}.
Then, if we consider a beginning of
spectral decomposition of $ K$   we have $
L^2\( \R , \C^2\) =\mathcal{N}\oplus \ker K  \oplus
\mathcal{P}$ with $\dim \mathcal{N}=1$ and $B>0$ in  $\mathcal{P}$ .
It is elementary to see that if $X$ is a vector subspace with     $\dim _\C X>3$ we have $X\cap
\mathcal{P}\neq 0$.   Suppose by contradiction that  $\ker \(  \mathcal{L} -\im \lambda   \)  \subsetneqq N_g\( \mathcal{L} -\im \lambda  \)$. Then there exist nontrivial  $u\in \ker \(  \mathcal{L} -\im \lambda   \)$ and $v\in N_g\( \mathcal{L} -\im \lambda  \)$ with $u= \(  \mathcal{L} -\im \lambda   \) ^nv$ for some $n\ge 1$. Then (by well known commutation properties  of $ \mathcal{L}$)
\begin{align*}
   \< Ku, u\> = \im \lambda  \< J^{-1}  u, u\> =  \im \lambda  \< J^{-1}  \(  \mathcal{L} -\im \lambda   \) ^nv, u\> = \im \lambda  \< J^{-1}v,   \(  \mathcal{L} -\im \lambda   \) ^nu\> =0.
\end{align*}
 So $B(u)=0$. Notice that we have also $B(\overline{u})=0$  where for the complex conjugate we have $\overline{u}\in \ker \(  \mathcal{L} +\im \lambda   \)$. Now consider the vector space
\begin{align*}
  X=\Span _\C \{ u , \overline{u} \} \oplus \ker K .
\end{align*}
We have $\dim _\C X =4$ with $B\le 0$ in $X$, which is a contradiction.

\qed

\section*{Acknowledgments}
C. was supported   by the Prin 2020 project \textit{Hamiltonian and Dispersive PDEs} N. 2020XB3EFL.
M.  was supported by the JSPS KAKENHI Grant Number 19K03579, 23H01079 and 24K06792.

Department of Mathematics, Informatics and Geosciences,  University
of Trieste, via Valerio  12/1  Trieste, 34127  Italy.
{\it E-mail Address}: {\tt scuccagna@units.it}

Department of Mathematics and Informatics,
Graduate School of Science,
Chiba University,
Chiba 263-8522, Japan.
{\it E-mail Address}: {\tt maeda@math.s.chiba-u.ac.jp}

\end{document}